\documentclass[12pt]{amsart}
\usepackage[mathscr]{eucal}
\usepackage{amsmath,amssymb,hyperref,lscape}
\usepackage{graphics}

\newtheorem{thm}{Theorem}[section]
\newtheorem{cor}[thm]{Corollary}
\newtheorem{prop}[thm]{Proposition}
\newtheorem{lem}[thm]{Lemma}

\newtheorem{remark}[thm]{Remark}
\newtheorem{example}[thm]{Example}

\begin{document}

\title[Godbillon-Vey sequences]
{Complex codimension one singular foliations and Godbillon-Vey sequences}
\author[CERVEAU, LINS-NETO, LORAY, PEREIRA and TOUZET]
{Dominique CERVEAU$^1$, Alcides LINS-NETO$^2$, Frank LORAY$^1$,\break 
Jorge Vit\'orio PEREIRA$^2$ and Fr\'ed\'eric TOUZET$^1$}
\address{\newline $1$ IRMAR, Campus de Beaulieu, 35042 Rennes Cedex, France\hfill\break
$2$ IMPA, Estrada Dona Castorina, 110, Horto, Rio de Janeiro, Brasil}
\email{dominique.cerveau,frank.loray,frederic.touzet@univ-rennes1.fr 
alcides,jvp@impa.br}
\date{September 2005}
\subjclass{}
\thanks{\newline
$1$ IRMAR, Campus de Beaulieu, 35042 Rennes Cedex, France\hfill\break
$2$ IMPA, Estrada Dona Castorina, 110, Horto, Rio de Janeiro, Brasil}
\keywords{}

\maketitle

\begin{abstract}
Let $\mathcal F$ be a codimension one singular holomorphic foliation
on a compact complex manifold $M$. Assume that there exists 
a meromorphic vector field $X$ on $M$ generically transversal to $\mathcal F$.
Then, we prove that $\mathcal F$ is the meromorphic pull-back
of an algebraic foliation on an algebraic manifold $N$,
or $\mathcal F$ is transversely projective outside a compact hypersurface, 
improving our previous work \cite{Crocodile}. 

Such a vector field insures the existence of a global meromorphic
Godbillon-Vey sequence for the foliation $\mathcal F$.
We derive sufficient conditions on this sequence 
insuring such alternative. For instance, if there exists a finite Godbillon-Vey
sequence or if the Godbillon-Vey invariant is zero, then either
$\mathcal F$ is the pull-back of a foliation on a surface,
or $\mathcal F$ is transversely projective. 
We illustrate these results with many examples.
\end{abstract}

\eject

\tableofcontents

%\eject

\section{Introduction}\label{S:intro}

Let $M$ be a compact connected complex manifold of dimension $n\ge2$. 
A (codimension $1$ singular holomorphic) foliation $\mathcal F$ on $M$ will be given 
by a covering of $M$ by open subsets $(U_j)_{j\in J}$ and a collection 
of integrable holomorphic $1$-forms $\omega_j$ on $U_j$, $\omega_j\wedge d\omega_j=0$,
having codimension $\ge2$ zero-set such that, on each non empty intersection $U_j\cap U_k$, we have
$$(*)\ \ \ \omega_j=g_{jk}\cdot\omega_k,\ \ \ \text{with}\ \ \ g_{jk}\in\mathcal
O^*(U_j\cap U_k).$$
Let $\text{Sing}(\omega_j)=\{p\in U_j\ ;\ \omega_j(p)=0\}$. Condition $(*)$
implies that $\text{Sing}(\mathcal F):=\cup_{j\in J}\text{Sing}(\omega_j)$
is a codimension $\ge2$ analytic subset of $M$. 
If $\omega$ is an integrable meromorphic $1$-form on $M$, $\omega\wedge d\omega=0$,
then we can associate to $\omega$ a foliation $\mathcal F_\omega$ as above.
Indeed, at the neighborhood of any point $p\in M$, one can write
$\omega=f\cdot\tilde{\omega}$ with $f$ meromorphic, sharing the same divisor 
with $\omega$; therefore, $\tilde{\omega}$ is holomorphic with codimension $\ge2$
zero-set and defines $\mathcal F_\omega$ on the neighborhood of $p$.

\eject

A {\it Godbillon-Vey sequence for $\mathcal{F}$}
is a sequence $(\omega_0,\omega_1,\ldots,\omega_k,\ldots)$
of meromorphic $1$-forms on $M$ 
such that $\mathcal{F}=\mathcal F_{\omega_0}$ and the formal $1$-form
\begin{equation}\label{E:integrability}
\Omega = dz + \sum_{k=0}^{\infty}  \frac{z^k}{k!}  \omega_k \, ,
\end{equation}
is integrable: $\Omega\wedge d\Omega=0$. In this sense, $\Omega$ defines
a formal development of $\mathcal{F}$ on the space $(\widehat{\mathbb{C},0})\times M$.
This condition is equivalent to
\begin{equation}\label{E:condition}
d \omega_k = \omega_0 \wedge \omega_{k+1} + 
\sum_{l=1}^k \binom{l}{k} \omega_l \wedge \omega_{k+1-l} \, .
\end{equation}
One can see that $\omega_{k+1}$ is well defined by $\omega_0,\omega_1,\ldots,\omega_k$
up to the addition by a meromorphic factor of $\omega_0$.
Conversally, $\omega_0,\omega_1,\ldots,\omega_k,\omega_{k+1}+f\cdot\omega_0$ 
is the begining
of another Godbillon-Vey sequence for any $f\in\mathcal M(M)$.

Given a meromorphic vector field $X$ on $M$ which is transversal to 
$\mathcal F$ at a generic point, there is a unique meromorphic $1$-form $\omega$
satisfying $\omega(X)=1$ and defining the foliation $\mathcal F$.
We then define a Godbillon-Vey sequence for $\mathcal F$ by setting
\begin{equation}\label{E:1}
\omega_k := L^{(k)}_X \omega,
\end{equation}
where $ L^{(k)}_X \omega$ denotes the $k^{\text{th}}$ Lie derivative along $X$ 
of the form $\omega$.
% (we have $d \Omega = \Omega \wedge  \frac{\partial \Omega}{\partial z}$).

The length of a Godbillon-Vey sequence is the minimal $N\in\mathbb N^*\cup\{\infty\}$ 
such that $\omega_k=0$ for $k\ge N$; in general, the length is infinite.
We say that $\mathcal F$ is {\it transversely projective}
if it admits a Godbillon-Vey sequence of length $\le3$, 
i.e. there are meromorphic $1$-forms $\omega_0=\omega$, $\omega_1$ and $\omega_2$ on $M$
satisfying
\begin{equation}\label{E:0}
\displaystyle{\left\lbrace
\begin{matrix}
d\omega_0 &=& \omega_0 \wedge \omega_1 \, \\
d\omega_1 &=& \omega_0 \wedge \omega_2 \, \\
d\omega_2 &=& \omega_1 \wedge \omega_2 \,
\end{matrix} \right. }
\end{equation}
This means that, outside the polar and singular set of the $\omega_i$'s,
the foliation $\mathcal F$ is (regular and) transversely projective in the classical sense
(see \cite{Godbillon} or section \ref{S:projectifRegular}) and this projective structure has 
``reasonable singularities''.
When $\omega_2=0$ (i.e. $d\omega_1=0$) or $\omega_1=0$ (i.e. $d\omega_0=0$), 
we respectively say that $\mathcal F_{\omega}$ is actually 
{\it transversely affine} or {\it euclidian}. 

\eject

Let $(\omega_k)$ be a Godbillon-Vey sequence for $\mathcal F$
and let $n$ be the smallest integer such that 
$\omega_0\wedge\cdots\wedge\omega_n\equiv0$, $1\le n\le m=\dim(M)$.
Then, the non trivial $n$-form $\Theta=\omega_0\wedge\cdots\wedge\omega_{n-1}$ is closed
and defines a singular codimension $n$ foliation $\mathcal F_\Theta$ 
whose leaves are contained in those of $\mathcal F$, $\mathcal F_\Theta\subset\mathcal F$.
We note that $\Theta$ does not depend on the choice of $\omega_1,\ldots,\omega_{n-1}$
in the Godbillon-Vey sequence, but does depend on  $\omega_0$.
Let $\mathcal M(M)$ be the field of meromorphic functions on $M$ and let
$K\subset\mathcal M(M)$ be the subfield of first integrals for $\mathcal F_\Theta$:
$$K=\{f\in\mathcal M(M)\ ;\ df\wedge\Theta\equiv0\}.$$
This field $K$ is integrally closed and, by \cite{Siegel}, there exists 
a meromorphic map $\pi:M\dashrightarrow N$ onto an algebraic manifold $N$
such that $K=\pi^*\mathcal M(N)$; in particular, the dimension $\dim(N)$ 
equals the transcendance degree of $K/\mathbb C$ and we have $1\le\dim(N)\le n\le m=\dim(M)$.
In the case $\Theta$ is a meromorphic volume form, that is $n=m$, 
we have $K=\mathcal M(M)$ and $N$ is the Algebraic Reduction of $M$ (see \cite{Ueno}).
We note that the fibration $\mathcal G$ induced on $M$ by the reduction map 
$\pi:M\dashrightarrow N$ contains $\mathcal F_\Theta$ as a sub-foliation
and may have any co\-dimension $n\le\dim(\mathcal G)\le a(M)\le m$, 
the algebraic dimension of $M$.
Our main theorem is the

\begin{thm}\label{T:Main}
Let $\mathcal F$ be a codimension $1$ singular foliation
on a compact complex manifold. Assume that $\mathcal F$
admits a global meromorphic Godbillon-Vey sequence $(\omega_k)$
and let $\Theta$, $K$ and $\pi:M\dashrightarrow N$ like above.
Then we are in one of the (non exclusive) following cases:
\begin{itemize}
\item $\mathcal F$ is the pull-back by $\pi:M\dashrightarrow N$ of a foliation $\underline{\mathcal F}$ on $N$,
\item or $\mathcal F$ is transversely projective.
\end{itemize}
\end{thm}

We are in the former case when the fibers of 
$\pi$ are contained in the leaves of $\mathcal F$;
this so happens for a generic foliation $\mathcal F$ on $M=\mathbb C\mathbb P(2)$:
$\pi$ is just the identity in this case.
Our statement becomes non trivial as soon as $M$ has not maximal 
algebraic dimension or when $\omega_0\wedge\cdots\wedge\omega_{m-1}\equiv0$,
$m=\dim(M)$.

When $N$ has dimension $n=0$ ($K=\mathbb C$) or $1$, 
then $\mathcal F$ is automatically transversely projective:
even in the case $\dim(N)=1$, the foliation $\underline{\mathcal F}$
has dimension $0$ and is trivially transversely euclidean.

\eject

We immediately deduce from Theorem \ref{T:Main} the

\begin{cor}\label{C:Main}
Let $\mathcal F$ be a codimension $1$ singular foliation
on a compact complex manifold. Assume that there exists
a meromorphic vector field $X$ on $M$ which is transversal to $\mathcal F$
at a generic point.
Then
\begin{itemize}
\item $\mathcal F$ is the pull-back by the Algebraic Reduction map
$M\dashrightarrow \text{red}(M)$ of a foliation on $N=\text{red}(M)$,
\item or $\mathcal F$ is transversely projective.
\end{itemize}
\end{cor}

In our previous work \cite{Crocodile}, this corollary was obtained 
under the stronger assumption
that the manifold $M$ is pseudo-parallelizable, i.e. there exist $m$
meromorphic vector fields $X_1,\ldots,X_m$ on $M$, $m=\dim(M)$,
that are independant at a generic point.

When $N$ has dimension $m-1$ or $m-2$, we prove that $\mathcal F$
is actually transversely affine if it is not a pull-back. In particular,
we have

\begin{thm}\label{T:GV=0}
Let $\mathcal F$ be a foliation on a compact complex manifold $M$ and 
assume that the meromorphic $3$-form $\omega_0 \wedge \omega_1 \wedge \omega_2$ is zero
for some Godbillon-Vey sequence associated to $\mathcal F$. Then
\begin{itemize}
\item $\mathcal F$ is the pull-back by a meromorphic map $\pi:M\dashrightarrow S$
of a foliation $\underline{\mathcal F}$ on an algebraic surface $S$,
\item or $\mathcal F$ is transversely affine.
\end{itemize}
\end{thm}

We do not know how to interpret this assumption geometrically. 
It is a well known fact and easy computation (see \cite{Godbillon}) 
that the meromorphic $3$-form $\omega_0 \wedge \omega_1 \wedge \omega_2$
is closed and well defined by $\mathcal F$ up to the addition 
by an exact meromorphic $3$-form. Nevertheless, we note that 
%suspensions of Hilbert modular 
foliations constructed in section \ref{S:modular}
%provide examples of transversely projective foliations with
have exact $3$-form $\omega_0 \wedge \omega_1 \wedge \omega_2$
but do not satisfy conclusion of Theorem \ref{T:GV=0}. 

Let us now define the length of a foliation, 
$\text{length}(\mathcal F)\in \mathbb{N^*} \cup \{ \infty \}$ 
as the minimal length among all Godbillon-Vey sequences attached to $\mathcal F$;
we set $\text{length}(\mathcal F)=\infty$ when $\mathcal F$ does not admit any Godbillon-Vey sequence.
A foliation has length $1$, $2$ or $3$ if, and only if, it is respectively 
transversely euclidian, affine or projective in the meromorphic sense above.
Also, consider an ordinary differential equation over a curve $C$
\begin{equation}\label{E:ODE}
  dz + \sum_{k=0}^{N}  \omega_k z^k\, ,
\end{equation}
(where $\omega_k$ are meromorphic $1$-forms defined on $C$).
Then, the foliation defined on $C \times \mathbb{C}P(1)$ by equation (\ref{E:ODE}) has
length $\le N+1$
(consider the Godbillon-Vey algorithm given by equation (\ref{E:1}) with $X={\partial\over\partial y}$).
Although it is expected that $N+1$ is the actual length of the generic equation (\ref{E:ODE}),
this is clear only for the Riccati equations ($N\le 2$), for monodromy reasons.

\eject

The study of foliations having finite length has been initiated by Camacho and Sc\'ardua 
in \cite{CamachoScardua} when the ambient space is a rational algebraic manifold.
We generalize their main result in the

\begin{thm}\label{T:CS}
Let  $\mathcal F$ be a foliation on a compact complex manifold $M$.
If $ 4\le \text{length}(\mathcal F) < \infty$, then $\mathcal F$ 
is the pull-back by a meromorphic map 
$\pi:M \dashrightarrow C \times \mathbb{C}P(1)$ of the foliation $\underline{\mathcal F}$ defined 
by an ordinary differential equation over a curve $C$ like above.
\end{thm}

There are examples of foliations on $\mathbb{C}P(2)$ having length $0$, $1$ or $2$
that are not pull-back of a Riccati equation (see \cite{LinsNeto} and \cite{Touzet}).
Therefore, condition $ 4\le \text{length}(\mathcal F)$ is necessary.
Recall that the degree of a foliation $\mathcal F$ on $\mathbb{C}P(n)$ 
is the number $d$ of tangencies with a generic projective line.
At least, we prove the

\begin{thm}\label{T:grau2}
Every foliation of degree $2$ on the complex projective space $\mathbb{C}P(n)$
has length at most $4$. This bound is sharp.
\end{thm}

In particular, Jouanolou examples (see \cite{Jouanolou}) have actually length $4$.
In the same spirit, we also derive from \cite{Loray} the

\begin{thm}\label{T:grau2loc}
If  $\mathcal F$ is a germ of foliation at the origin of $\mathbb{C}^n$ defined 
by an holomorphic $1$-form with a non zero linear part, then $\text{length}(\mathcal F) \le 4$.
\end{thm}

From Theorems \ref{T:CS} and \ref{T:grau2}, we immediately retrieve 
the following result previously obtained by two of us in \cite{CerveauLinsNeto}: 

\begin{cor}\label{C:CerveauLinsNeto} A degree $2$ foliation
on $\mathbb{C}P(n)$ is either transversely projective, or the pull-back of a foliation
on $\mathbb{C}P(2)$ by a rational map.\end{cor}

We do not understand the strength of the assumption $\text{length}(\mathcal F) < \infty$ of Theorem
\ref{T:CS}. In fact,
we still do not know any example of a foliation having finite length $>4$.
It is not excluded that the generic foliation of degree $3$ on $\mathbb{C}P(2)$
has infinite length.

In section \ref{S:Examples}, we also provide examples of transversely
projective foliations on $\mathbb{C}P(3)$ that are not transversely
affine. In fact, they form a new irreducible component of the space
of foliations of degree $4$ (see \cite{CerveauLinsNeto}).
We do not know yet if they are pull-back by rational map
of foliations on $\mathbb CP(2)$.
We also give an example of a degree $6$ transversely
projective foliation $\mathcal H_2$ in $\mathbb{C}P(3)$ (with explicit equations) 
which is not the pull-back of a foliation in $\mathbb{C}P(2)$ by a rational map. 
In fact, $\mathcal H_2$ is the suspension (see section \ref{S:projectifSingular}) 
of one of the ``Hilbert modular foliations'' on $\mathbb{C}P(2)$ studied 
in \cite{MendesPereira}. We do not know if this foliation is isolated 
in the space of foliations.

\eject

Finally, since our arguments are mainly of algebraic nature, it is natural to ask
what remains true from our work in the positive characteristic. In this direction,
we prove in the last section the 

\begin{thm}\label{T:positiva}
Let $M$ be a smooth projective variety defined over a field $K$ 
of characteristic $p > 0$ and $\omega$ be a rational $1$-form. 
If $\omega$ is integrable $\omega\wedge d\omega=0$, then there exist a rational function 
$F \in K(M)$ such that $F\omega$ is closed. In this sense, 
the "foliation" $\mathcal F_\omega$ has length $1$.
\end{thm}

%%%%%%%%%%%%%%%%%%%%%%%%%%%%%%%%%%%%%%%%%%%%%%%%%%%%%%%%%%%%%%%%%%%%%%%%%%%%%%%%%%%%%%%%%%%%%%%%%%%%%%%%%%%%%%%%%%%%%%%%%%%%%%%%

\section{Background and first steps}\label{S:Material}

\subsection{Godbillon-Vey sequences \cite{Godbillon,CamachoScardua}}\label{S:GodbillonVey}

We introduce Godbillon-Vey sequences 
for a codimension one foliation $\mathcal F$
and describe basic properties.
Let $\omega$ be a differential $1$-form defining $\mathcal F$ and $X$ be a vector field
satisfying $\omega(X)=1$. Then, the integrability condition of $\omega$
is equivalent to
\begin{equation}\label{E:int}
\omega\wedge d\omega=0\ \ \ \Leftrightarrow\ \ \ d\omega=\omega\wedge L_X\omega.
\end{equation}
Indeed, from $L_X\omega=d(\omega(X))+d\omega(X,.)=d\omega(X,.)$, we derive
$$0=\omega\wedge d\omega(X,.,.)=\omega(X)\cdot d\omega-\omega\wedge(d\omega(X,.))
=d\omega-\omega\wedge L_X\omega$$
(the converse is obvious). Applying this identity to the formal $1$-form
\begin{equation}\label{E:DefinitionOmega}
\Omega=dz+ \omega_0+z\omega_1+{z^2\over2}\omega_2+\cdots+\frac{z^k}{k!}\omega_k+\cdots
\end{equation}
together with the vector field $X=\partial_z$, we derive
$$\Omega\wedge d\Omega=0\ \ \ \Leftrightarrow\ \ \ 
\sum_{k=0}^{\infty}\frac{z^k}{k!}  d\omega_k
=\left(\sum_{k=0}^{\infty}\frac{z^k}{k!}\omega_k\right)\wedge
\left(\sum_{k=1}^{\infty}\frac{z^{k-1}}{(k-1)!}\omega_k\right).$$
We therefore obtain the full integrability condition (\ref{E:condition}) for $\Omega$:
\begin{equation}\label{E:FullIntegrability}
\begin{matrix}
d\omega_0&=&\omega_0\wedge\omega_1\hfill\\
d\omega_1&=&\omega_0\wedge\omega_2\hfill\\
d\omega_2&=&\omega_0\wedge\omega_3+\hskip0.2cm\omega_1\wedge\omega_2\hfill\\
d\omega_3&=&\omega_0\wedge\omega_4+2\omega_1\wedge\omega_3\hfill\\
&\vdots&\\
d\omega_k&=& \omega_0 \wedge \omega_{k+1} + 
\sum_{l=1}^k \binom{l}{k} \omega_l \wedge \omega_{k+1-l} \\
&\vdots&
\end{matrix}
\end{equation}
For instance, if we start with $\omega$ integrable and $X$ satisfying $\omega(X)=1$,
then the iterated Lie derivatives $\omega_k := L^{(k)}_X \omega$ 
define a Godbillon-Vey sequence for $\mathcal F_{\omega}$.
Indeed, from the formula $(L_X\omega)(X)=d\omega(X,X)=0$, we have
$$\omega_0(X)=1\ \ \ \text{and}\ \ \ \omega_k(X)=0\ \text{for all}\ k>0;$$
therefore, $\Omega(X)=1$ and integrability condition comes from
$$\Omega\wedge L_X\Omega
=\left(dz+\sum_{k=0}^{\infty}\frac{z^k}{k!}\omega_k\right)\wedge
\left(\sum_{k=0}^{\infty}\frac{z^k}{k!}\omega_{k+1}\right)
=d\Omega.$$

%\eject

From a given Godbillon-Vey sequence, we derive many other ones. 
For instance, given any non zero 
meromorphic function $f\in\mathcal M(M)$, after applying 
the formal change of variable $z=f\cdot t$ to 
$$\Omega=dz+ \omega_0+z\omega_1+{z^2\over2}\omega_2+\cdots+\frac{z^k}{k!}\omega_k+\cdots,$$
we derive the new integrable $1$-form
\[
  {\Omega\over f} = dt +  \frac{\omega_0}{f} + t( \omega_1 + \frac{df}{f}) + 
  \frac{t^2 }{2}(f\omega_2) + \frac{t^3 }{3!}(f^2\omega_3) +\cdots +
  \frac{t^k }{k!}(f^{k-1}\omega_k)+\cdots
\]
In other words, we obtain a new Godbillon-Vey sequence $(\tilde\omega_k)$ by setting
\begin{equation}\label{E:changeOmega0}
\left\{\begin{matrix}
\tilde\omega_0&=&{1\over f}\cdot\omega_0\hfill\\
\tilde\omega_1&=&\hskip0.5cm\omega_1 + \frac{df}{f}\\
\tilde\omega_2&=&f\cdot\omega_2\hfill\\
&\vdots&\\
\tilde\omega_{k+1}&=&f^k\cdot\omega_{k+1}\hfill\\
&\vdots&
\end{matrix}\right.
\end{equation}
By the same way, we can apply to $\Omega$ the formal change of variable 
$z=t+f\cdot t^{k+1}$, 
$k=1,2,\ldots$, and successively derive new Godbillon-Vey sequences
\begin{equation}\label{E:changeOmegak}
\left\{\begin{matrix}
\tilde\omega_0&=&\omega_0\hfill\\
\tilde\omega_1&=&\omega_1 + f\omega_0\hfill\\
\tilde\omega_2&=&\omega_2+f\omega_1-df\hfill\\
&\vdots&\\
\end{matrix}\right.\ \ \ \left\{\begin{matrix}
\tilde\omega_0&=&\omega_0\hfill\\
\tilde\omega_1&=&\omega_1 \hfill\\
\tilde\omega_2&=&\omega_2+f\omega_0\hfill\\
&\vdots&\\
\end{matrix}\right.\ \ \ \text{etc\ldots}
\end{equation}
Conversally, we easily see from integrability condition (\ref{E:condition})
that $\omega_{k+1}$ is well defined by $\omega_0,\omega_1,\ldots,\omega_k$
up to the addition by a meromorphic factor of $\omega_0$.
In fact, every Godbillon-Vey sequence can be deduced from a given one 
after applying to the $1$-form $\Omega$ a formal transformation 
belonging to the following group
\[
  G = \left\{ (p,z) \mapsto \left(p, \sum_{k=1}^{\infty} f_k(p)\cdot z^k\right),\ 
  f_k\in\mathcal M(M),\  f_1 \not\equiv 0 \right\}.
\]
In particular, the so-called Godbillon-Vey invariant
$\omega_0\wedge\omega_1\wedge\omega_2=-\omega_1\wedge d\omega_1$ 
is closed and is well defined up to the addition 
by an exact meromorphic $3$-form of the form
$${df\over f}\wedge\omega_0\wedge\omega_2={df\over f}\wedge d\omega_1 \ \ 
\text{or}\ \ \ 
df\wedge\omega_0\wedge\omega_1=df\wedge d\omega_0$$
for some meromorphic function $f\in\mathcal M(M)$.

\begin{remark}\rm
A natural Godbillon-Vey sequence for the formal foliation 
$\mathcal F_\Omega$ defined by $\Omega$ is given by 
$$\Omega_k=L_{\partial_z}^{(k)}\Omega=\sum_{l=k}^{\infty}\frac{z^{l-k}}{(l-k)!}\omega_l,\ \ \ 
k>0$$
or equivalently by the formal integrable $1$-form
$$d(t+z)+\omega_0+(t+z)\omega_1+{(t+z)^2\over2}\omega_2+\cdots$$
$$=dt+\Omega_0+t\Omega_1+{t^2\over2}\Omega_2+\cdots$$
\end{remark}

In fact, this remark also applies to the case where the $\omega_k$
are meromorphic $1$-forms on a complex curve $C$.
The so-called ``ordinary differential equation'' defined by
$$
  \Omega=dz + \sum_{k=0}^{N} {z^k\over k!} \omega_k \, ,
$$
defines a foliation $\mathcal F$ on $C \times \mathbb{C}P(1)$
(integrability conditions (\ref{E:condition}) are trivial in dimension $1$).
This foliation admits a natural Godbillon-Vey sequence of length $N+1$
given by $L_{\partial_z}^{(k)}\Omega$ (or by replacing $z$ by $z+t$).

\begin{remark}\rm
It follows from relations (\ref{E:condition}) that all differential forms 
$$
\Theta_k:=\omega_0\wedge\omega_1\wedge\cdots\wedge\omega_{k-1}\ \ \ 
\text{for all}\ k=2,\ldots,n,
$$ 
are closed and depend only on $\omega_0$. We obtain an ``integrable flag'':
$$
\mathcal F=\mathcal F_0\supset\mathcal F_1\supset\cdots\supset\mathcal F_{n-1}
$$
(the tangents spaces $T_p\mathcal F_k$ define is a flag at a generic point $p\in M$).
The codimension $n$ of the flag is the first $n$ such that 
$\omega_0\wedge\cdots\wedge\omega_n=0$.
\end{remark}

\eject

We have two preliminary lemmas about finite Godbillon-Vey sequences.

\begin{lem}\label{L:tangGVfini}
Let $\omega_0,\omega_1,\ldots,\omega_N$ 
be a Godbillon-Vey sequence of finite length $N+1$. 
Then $\omega_k\wedge\omega_l=0$ for all $k,l\ge2$ 
and integrability conditions become
$$d\omega_k=\omega_0\wedge\omega_{k+1}+(k-1)\omega_1\wedge\omega_k\ \ \ 
k=0,1,\ldots,N.$$
\end{lem}

In particular, the condition $\omega_{N+1}=0$ in a Godbillon-Vey
sequence is not sufficient to conclude that the truncated sequence
$$\omega_0,\omega_1,\ldots,\omega_N,0,0,\ldots$$
provides a finite Godbillon-Vey sequence, except when $N=0,1$ or $2$.

\begin{proof}We assume $\omega_N\not=0$ with $N\ge2$, otherwise we have done. 
The integrability conditions (\ref{E:condition}) 
$$\begin{matrix}
&d\omega_0&=&\omega_0\wedge\omega_1\hfill\\
&d\omega_1&=&\omega_0\wedge\omega_2\hfill\\
&&\vdots&\\
&d\omega_N&=& \sum_{l=1}^N \binom{l}{N} \omega_l \wedge \omega_{N+1-l} \\
0=&d\omega_{N+1}&=&\sum_{l=2}^N \binom{l}{N+1} \omega_l \wedge \omega_{N+2-l} \\
&&\vdots&\\
0=&d\omega_{2N-2}&=&{1\over N}\binom{N-1}{2N-2}\omega_{N-1}\wedge\omega_N\hfill
\end{matrix}$$
Examining the line of index $k=2N-2$,
we deduce that $\omega_{N-1}\wedge \omega_N \equiv 0$. 
Futhermore, by descendent induction, we also deduce from the line of index $k+N-1$
that $\omega_{k}\wedge \omega_N \equiv 0$ for every $k\ge 2$. 
Therefore, the remining $N$ first lines of integrability conditions
are as in the statement.
\end{proof}

\begin{cor}\label{C:uniciteGVfini}
Let $\omega_0$, $\omega_1$ and $\omega_2$ be differential $1$-forms
satisfying relations (\ref{E:condition}) for $k=0,1$ with $d\omega_1\neq0$.
Then, there exists at most one finite Godbillon-Vey sequence 
$\omega_0,\ldots,\omega_N$ completing this triple.
\end{cor}

\begin{proof}The assumption $d\omega_1=\omega_0\wedge\omega_2\neq0$ 
implies in particular that $\omega_2\neq0$.
 If $\omega_0,\omega_1,\ldots,\omega_N$ is a finite sequence, 
then we recursively see from integrability conditions of Lemma \ref{L:tangGVfini}
that the line of index $k$ determines $\omega_k$, $k=3,\ldots,N$, 
up to a meromorphic factor of $\omega_0$. But since $\omega_k$ is tangent 
to $\omega_2$ but $\omega_0$ is not, we deduce that $\omega_k$ is actually
completely determined by the line of index $k$.
\end{proof}

\eject

Here is a weaker but easier version of Theorem \ref{T:CS}.

\begin{thm}\label{T:weakCS}
Let $\mathcal F$ be a foliation on a compact pseudo-parallelizable manifold $M$. 
If $\text{length}(\mathcal F)<\infty$, then we have the following alternative:
\begin{enumerate}
\item $\mathcal F$ is the pull-back
of a foliation $\underline{\mathcal F}$ on an algebraic surface $S$
by a meromorphic map $\pi: M \dashrightarrow S$ with 
$\text{length}(\underline{\mathcal F})=\text{length}(\mathcal F)$,
\item or $\mathcal F$ is transversely projective, 
i.e. $\text{length}(\mathcal F)\le 3$.
\end{enumerate}
\end{thm}

\begin{proof}Let $(\omega_0, \omega_1, \ldots, \omega_N)$ be a Godbillon-Vey 
sequence for $\mathcal F$ with $\omega_N\not=0$, $N\ge3$
and $\omega_1$, $\omega_2$, $\omega_3$ both non zero
(otherwise we are in the second alternative of the statement). 
Following Lemma \ref{L:tangGVfini}, there exist meromorphic functions $f_k$ 
such that $\omega_k=f_k \cdot \omega_2$. Observe that $f_3 \ne 0$ since $\omega_3\ne 0$.
Recall that $( \omega_k )$ is a Godbillon-Vey sequence and the $1$-form
\[
 \Omega = dz + \omega_0 + z \omega_1 + \cdots + \frac{z^N}{N!}\omega_N \, ,
\]
is integrable. Applying to $\Omega$ the change of variables $z=t/f_3$
(see Section \ref{S:GodbillonVey}), we derive a new Godbillon-Vey sequence
of length $N$ satisfying $\omega_2=\omega_3$. Therefore
\[
\displaystyle{\left\lbrace
\begin{matrix}
d\omega_2 &=& \omega_0 \wedge \omega_{3} &+& \hfill\omega_1 \wedge \omega_{2} 
&=&\hfill \omega_0 \wedge \omega_{2} &+&\hfill \omega_1 \wedge \omega_{2} \\
d\omega_3 &=& \omega_0 \wedge \omega_{4} &+& 2\cdot  \omega_1 \wedge \omega_{3}  
&=& f_4\cdot \omega_0 \wedge \omega_{2} &+& 2\cdot \omega_1 \wedge \omega_{2}
\end{matrix} \right. }
\]
In particular $(1-f_4)\omega_0 \wedge \omega_2 = \omega_1 \wedge \omega_2$ 
implying that $\omega_0\wedge \omega_1 \wedge \omega_2 \equiv 0$.
We conclude with Theorem \ref{T:GV=0} (a consequence of section \ref{S:main}).

Let us now prove that $\text{length}(\underline{\mathcal F})=\text{length}(\mathcal F)$.
Since a Godbillon-Vey sequence
for $\underline{\mathcal F}$ induces, by pull-back by $\phi$, a sequence for $\omega_0$, 
it follows that $\text{length}(\underline{\mathcal F}) \ge \text{length}(\mathcal F) = N$.
Let $\underline\omega_0$ be the meromorphic $1$-form on $S$ 
such that $\phi^*\underline\omega_0 = \omega_0$. 
From the equality $0=\omega_0 \wedge \omega_1 \wedge\omega_2 = \omega_1 \wedge  d\omega_1$, 
we see that $\omega_1$ is integrable.
Writing down the equations in local coordinates we also see that  the fibers of $\phi$
are tangent to  the  foliation associated to $\omega_1$. 
Moreover, $\omega_1$ is the pull-back by
$\phi$ of a $1$-form $\underline\omega_1$ on $S$.
Recall that $\omega_2 = f_0 \omega_0 + f_1 \omega_1$ 
and that  $df_1\wedge d\omega_0=0$.
Differentiating  the identity
\[
  d \omega_2 = \omega_0 \wedge \omega_2 + \omega_1 \wedge \omega_2 = (f_1 - f_0) d \omega_ 0
\]
it follows that  $df_0 \wedge d\omega_0=0$. 
Consequently $\omega_2 = \phi^* \underline\omega_2$, where  $\underline\omega_2$ is a meromorphic
$1$-form on $S$.
At this point we can rewrite $\Omega$ as
\[
\Omega = dz + \phi^* \underline\omega_0 + z \phi^* \underline\omega_1 + h \cdot \phi^* \underline\omega_2 \, ,
\]
where $h = \frac{z^2}{2} + \sum_{i=3}^N \frac{z^N}{N!} h_n$.
The integrability of $\Omega$ implies that $dh \wedge \phi^* \underline\omega_2 = 0$, 
where $d$ is the differential over $M$ (i.e. $dz=0$).
This implies that each $h_j$ belongs to $\phi^{-1}\mathcal M(S)$ 
and therefore $\omega_j = \phi^* \underline\omega_j$ for every $j$ and some $\underline\omega_j$
on $S$. This proves that $\text{length}(\underline{\mathcal F})\le N$.
\end{proof}

\eject

%%%%%%%%%%%%%%%%%%%%%%%%%%%%%%%%%%%%%%%%%%%%%%%%%%%%%%%%%%%%%%%%%%%%%%%%%%%%%%%%%%%%%%%%%%%%%%%%%%%%%%%%%%%%%%%%%%%%%%%%%%%%%%%%

\subsection{Transversely projective foliations: the classical case
\cite{Godbillon,SullivanThurston}}\label{S:projectifRegular}

A {\bf regular} codimension one foliation $\mathcal F$ on a manifold $M$ 
is transversely projective if there exists an atlas of submersions
$f_i:U_i \to \mathbb CP(1)$ on $M$ satisfying the cocycle condition:
\[
f_i =\frac{ a_{ij} f_j + b_{ij}}{c_{ij} f_j + d_{ij}}, \ \ \  
\left( \begin{matrix} a_{ij} & b_{ij} \\ c_{ij} & d_{ij} \end{matrix} \right)\in PGL(2,\mathbb C).
\]
on any intersection $U_i\cap U_j$. Any two such atlases $(f_i:U_i \to \mathbb CP(1))_i$
and $(g_k:V_k \to \mathbb CP(1))_k$ define the same projective structure 
if the union of them is again a projective structure, i.e. satisfying
the cocycle condition $f_i =\frac{ a_{ik} g_k + b_{ik}}{c_{ik} g_k + d_{ik}}$
on $U_i\cap V_k$. 

Starting from one of the local submersions $f:U\to \mathbb CP(1)$ above, 
one can step-by-step modify the other charts so that
they glue with $f$ and define an analytic continuation for $f$. 
Of course, doing this along an element $\gamma\in\pi_1(M)$ 
of the fundamental group, we obtain monodromy $f(\gamma\cdot p)=A_\gamma\cdot f(p)$
for some $A_\gamma\in  PGL(2,\mathbb C)$. By this way, we define 
the {\it monodromy representation} of the structure,
that is a homomorphism
$$\rho:\pi_1(M)\to PGL(2,\mathbb C);\gamma\mapsto A_\gamma,$$
as well as the {\it developing map}, 
that is the full analytic continuation of $f$ 
on the universal covering $\tilde M$ of $M$
$$\tilde f:\tilde M\to \mathbb CP(1).$$
By construction, $\tilde f$ is a global submersion on $\tilde M$
whose determinations $f_i:U_i \to \mathbb CP(1)$ on simply connected
subsets $U_i\subset M$ define unambiguously the foliation $\mathcal F$
and the projective structure. In fact, the map $\tilde f$
is $\rho$-equivariant 
\begin{equation}\label{E:rhoEquiv}
f(\gamma\cdot p)=\rho(\gamma)\cdot f(p),\ \ \ \forall\gamma\in\pi_1(M).
\end{equation}
Finally, we obtain

\begin{prop}A regular foliation $\mathcal F$ on $M$ is transversely
projective if, and only if, there exist
\begin{itemize}
\item a representation $\rho:\pi_1(M)\to PGL(2,\mathbb C)$ 
\item a submersion $\tilde f:\tilde M\to \mathbb CP(1)$ defining $\mathcal F$ 
and satisfying (\ref{E:rhoEquiv}).
\end{itemize} 
Any other pair $(\rho',\tilde f')$
will define the same structure if, and only if, we have
$\rho'(\gamma)=A\cdot\rho(\gamma)\cdot A^{-1}$ and $\tilde f'=A\cdot\tilde f$.
for some $A\in PGL(2,\mathbb C)$.\end{prop}

\eject

\begin{remark}\label{R:simplyConnected}\rm
If $M$ is simply connected, then any transversely projective foliation
$\mathcal F$ on $M$ actually admits a global first integral
$\tilde f:M\to \mathbb CP(1)$, a holomorphic mapping.
\end{remark}

\begin{example}[\it Suspension of a representation]\rm
Given a representation $\rho:\pi_1(M)\to PGL(2,\mathbb C)$ 
of the fundamental group of a manifold $M$ into the projective group,
we derive the following representation into the group of diffeomorphisms
of the product $\tilde M\times \mathbb CP(1)$
$$\tilde\rho:\pi_1(M)\to\text{Aut}(\tilde M\times \mathbb CP(1))\ ;\ 
(p,z)\mapsto(\gamma\cdot p,\rho(\gamma)\cdot z)$$
($\tilde M$ is the universal covering of $M$ and $p\mapsto\gamma\cdot p$, 
the Galois action of $\gamma\in\pi_1(M)$). 
The image $\tilde G$ of this representation acts freely, properly and discontinuously
on the product $\tilde M\times \mathbb CP(1)$ since its restriction to the first factor does.
Moreover, $\tilde G$ preserves the horizontal foliation $\mathcal H$
defined by $dz$ as well as the vertical $\mathbb CP(1)$-fibration defined
by the projection $\pi:\tilde M\times \mathbb CP(1)\to\tilde M$ onto the first factor.
In fact, we have 
$\pi(\tilde\rho(\gamma)\cdot p)=\rho(\gamma)\cdot \tilde\pi(p)$
for all $p\in\tilde M$ and $\gamma\in\pi_1(M)$.
Therefore, the quotient $N:=\tilde M\times \mathbb CP(1)/\tilde G$ is a manifold
equipped with a locally trivial $\mathbb CP(1)$-fibration given by the projection 
$\pi:N\to M$ as well as a codimension one foliation $\mathcal H$ transversal
to $\pi$. In fact, the foliation $\mathcal H$ is transversely projective
with monodromy representation $\rho\circ\pi_*:\pi_1(N)\to PGL(2,\mathbb C)$ 
($\pi$ induces an isomorphism $\pi_*:\pi_1(N)\to\pi_1(M)$) 
and developing map $\tilde M\times \mathbb CP(1)\to \mathbb CP(1);(p,z)\mapsto z$
($\tilde M\times \mathbb CP(1)$ is the universal covering of $N$).

Conversely, a codimension one foliation $\mathcal H$
transversal to a $\mathbb CP(1)$-fibration $\pi:N\to M$ 
is actually the suspension
of a representation $\rho:\pi_1(M)\to PGL(2,\mathbb C)$.
In particular, $\mathcal H$ is transversely projective 
and uniquely defined by its monodromy $\rho$.
\end{example}

Now, given a transversely projective foliation $\mathcal F$ on $M$,
we construct {\it the suspension of $\mathcal F$} as follows.
We first construct the suspension of the monodromy representation 
$\rho:\pi_1(M)\to PGL(2,\mathbb C)$ of $\mathcal F$ as above
and consider the graph
$$\tilde\Gamma=\{(p,z)\in\tilde M\times \mathbb CP(1)\ ;\ z=\tilde f(p)\}$$
of the developing map $\tilde f:\tilde M\to \mathbb CP(1)$. 
Since $\tilde f$ is $\rho$-equivariant, 
its graph $\tilde\Gamma$ is invariant under the group $\tilde G$ and defines
a smooth cross-section $f:M\hookrightarrow N$ to the $\mathbb CP(1)$-fibration
$\pi:N\to M$. By construction, its image $\Gamma=f(M)$ is also
transversal to the ``horizontal foliation'' $\mathcal H$ 
and the transversely projective foliation induced by $\mathcal H$ on $\Gamma$ 
actually coincides (via $f$ or  $\pi$)
with the initial foliation $\mathcal F$ on $M$.

\eject

\begin{prop}\label{P:Suspension}
A regular foliation $\mathcal F$ on $M$ is transversely
projective if, and only if, there exist
\begin{itemize}
\item a locally trivial $\mathbb CP(1)$-fibration $\pi:N\to M$ over $M$,
\item a codimension one foliation $\mathcal H$ on $N$ transversal to $\pi$,
\item a section $f:M\to N$ transversal to $\mathcal H$ such that
the foliation induced by $\mathcal H$ on $f(M)$ coincides via $f$ with $\mathcal F$.
\end{itemize} 
Any other triple $(\pi':N'\to M,\mathcal H',f')$
will define the same structure if, and only if, there exists
a diffeomorphism $\Phi:N'\to N$ such that $\pi'=\pi\circ\Phi$,
$f=\Phi\circ f'$ and $\mathcal H'=\Phi^*\mathcal H$.\end{prop}

Over any sufficiently small open subset $U\subset M$, 
the $\mathbb CP(1)$-fibration is trivial and
one can choose trivializing coordinates $(p,z)\in U\times \mathbb CP(1)$ 
such that $f:U\to \pi^{-1}(U)$ coincides with the zero-section $\{z=0\}$.
The foliation $\mathcal H$ is defined by a unique differential $1$-form
of the type
$$\Omega=dz+\omega_0+z\omega_1+z^2\omega_2$$
where $\omega_0$, $\omega_1$ and $\omega_2$ are holomorphic $1$-forms 
defined on $U$. The integrability condition $\Omega\wedge d\Omega=0$
reads
\begin{equation}\label{E:triplet}
\left\{\begin{matrix}
d\omega_0 &=& \hfill\omega_0 \wedge \omega_1 \\
d\omega_1 &=& 2\omega_0 \wedge \omega_2 \\
d\omega_2 &=& \hfill\omega_1 \wedge \omega_2 
\end{matrix} \right. 
\end{equation}
Now, any change of trivializing coordinates preserving the zero-section
takes the form $(\tilde p,\tilde z)=(p,f_0\cdot z/(1+f_1\cdot z))$ 
where $f_0:U\to\mathbb C^*$ and $f_1:U\to\mathbb C$ are holomorphic.
The foliation $\mathcal H$ is therefore defined by
$$\tilde\Omega:={(f_0-f_1\tilde z)^2\over f_0}\Omega=
d\tilde z+\tilde \omega_0+\tilde z\tilde \omega_1+\tilde z^2\tilde \omega_2$$
where the new triple $(\tilde\omega_0,\tilde\omega_1,\tilde\omega_2)$ 
is given by
\begin{equation}\label{E:tripletModif}
\left\{\begin{matrix}
\tilde\omega_0 &=& f_0\omega_0 \hfill \\
\tilde\omega_1 &=& \omega_1 -2f_1\omega_0-{df_0\over f_0} \hfill\\
\tilde\omega_2 &=& {1\over f_0}\left(\omega_2-f_1\omega_1+f_1^2\omega_0 +df_1\right)
\end{matrix} \right.
\end{equation}

\begin{prop}A regular foliation $\mathcal F$ on $M$ is transversely
projective if, and only if, there exists an atlas of charts $U_i$ 
equipped with $1$-forms $(\omega_0^i,\omega_1^i,\omega_2^i)$ 
satisfying (\ref{E:triplet}) and related to each other
by (\ref{E:tripletModif}) on $U_i\cap U_j$.
\end{prop}

\eject

\begin{example}\label{E:SL2}\rm
Consider 
$$SL(2,\mathbb C)=\{\left(\begin{matrix} x&u\\y&v\end{matrix}\right)\ ;\ xv-yu=1\}.$$ 
The meromorphic function defined by
$$f:SL(2,\mathbb C)\to \mathbb CP(1)\ ;\ \left(\begin{matrix} x&u\\y&v\end{matrix}\right)
\mapsto {x\over y}$$
is a global submersion defining a transversely projective foliation $\mathcal F$
on $SL(2,\mathbb C)$. The leaves are the 
right cosets for the ``affine'' subgroup
$$\mathbb A=\{\left(\begin{matrix} a&b\\0&{1\over a}\end{matrix}\right)\ ;\ a\not=0\}.$$
Indeed, we have for any $z\in\mathbb C$
$$\left(\begin{matrix} z&-1\\1&0\end{matrix}\right)\cdot\mathbb A=
\{\left(\begin{matrix} az&bz-{1\over a}\\a&b\end{matrix}\right)\ ;\ a\not=0\}
=\{f=z\}$$
and for any $w=1/z\in\mathbb C$
$$\left(\begin{matrix} 1&0\\w&1\end{matrix}\right)\cdot\mathbb A=
\{\left(\begin{matrix} a&b\\aw&bw+{1\over a}\end{matrix}\right)\ ;\ a\not=0\}
=\{f=1/w\}.$$
In fact, if we consider the projective action
of a matrix $\left(\begin{matrix} x&u\\y&v\end{matrix}\right)$ 
on $(z:1)\in \mathbb CP(1)$, then $f$ is nothing but the image of the direction $(1:0)$
(i.e. $z=\infty$) by the matrix and $\{f=\infty\}$ coincides with the affine
subgroup $\mathbb A$ fixing $z=\infty$.

A global holomorphic triple $(\omega_0,\omega_1,\omega_2)$ for $\mathcal F$ 
can be constructed as follows. Consider the Maurer-Cartan form
$$\mathcal M:=\left(\begin{matrix} x&u\\y&v\end{matrix}\right)^{-1}
d\left(\begin{matrix} x&u\\y&v\end{matrix}\right)
=\left(\begin{matrix} vdx-udy&vdu-udv\\xdy-ydx&xdv-ydu\end{matrix}\right).$$
The matrix $\mathcal M$ is a differential $1$-form on $SL(2,\mathbb C)$ 
taking values in the Lie algebra $sl(2,\mathbb C)$ 
($\text{trace}(\mathcal M)=d(xv-yu)=0$)
and its coefficients form a basis for the left-invariant $1$-forms 
on $SL(2,\mathbb C)$. If we set
$$\mathcal M=\left(\begin{matrix} 
-{\omega_1\over 2}&-\omega_2\\ \omega_0&{\omega_1\over 2}
\end{matrix}\right),$$
then Maurer-Cartan formula
$d\mathcal M+\mathcal M\wedge\mathcal M=0$
is equivalent to integrability conditions (\ref{E:triplet}) for the triple
$(\omega_0,\omega_1,\omega_2)$. 
In fact, the ``meromorphic triple'' 
$$(\tilde\omega_0,\tilde\omega_1,\tilde\omega_2)=(df,0,0)$$
is derived by setting $f_0=-{1\over y^2}$ 
and $f_1=-{v\over y}$ in formula (\ref{E:tripletModif}).

\eject

A left-invariant $1$-form $\omega=\alpha\omega_0+\beta\omega_1+\gamma\omega_2$, 
$\alpha,\beta,\gamma\in\mathbb C$, is integrable, $\omega\wedge d\omega=0$, 
if, and only if, $\alpha\gamma=\beta^2$.
The right translations act transitively on the set of integrable 
left-invariant $1$-forms and thus on the corresponding foliations. 
For instance, if we denote by $T_z$ the right translation
$$T_z:SL(2,\mathbb C)\to SL(2,\mathbb C)\ ;\ 
\left(\begin{matrix} x&u\\y&v\end{matrix}\right)\mapsto
\left(\begin{matrix} x&u\\y&v\end{matrix}\right)
\left(\begin{matrix} z&-1\\1&0\end{matrix}\right),\ \ \ z\in\mathbb C,$$
then we have $T_z^*\omega_0=z^2\omega_0+z\omega_1+\omega_2$
and the corresponding foliation $\mathcal F_z$ is actually defined
by the global submersion 
$$f\circ T_z:SL(2,\mathbb C)\to \mathbb CP(1)\ ;\ 
\left(\begin{matrix} a&b\\c&d\end{matrix}\right)\mapsto
{az+b\over cz+d}.$$ 
The leaf $\{f\circ T_z=w\}$ of $\mathcal F_z$ is the set of matrices
sending the direction $(z:1)$ onto $(w:1)$.
\end{example}

\begin{remark}\rm
Let $(\omega_0,\omega_1,\omega_2)$ be a triple of holomorphic $1$-forms
on a manifold $M$ satisfying integrability condition (\ref{E:triplet}).
The differential equation
$$dz+\omega_0+z\omega_1+z^2\omega_2=0$$
defined on the trivial projective bundle $M\times \mathbb CP(1)$ can be lifted
as an integrable differential $sl(2,\mathbb C)$-system defined on 
the rank $2$ vector bundle $M\times\mathbb C^2$ by
$$\left\{\begin{matrix} dz_1&=&-{\omega_1\over 2}z_1&-&\omega_2z_2\\
dz_2&=&\hfill\omega_0z_1&+&\hfill{\omega_1\over 2}z_2\end{matrix}\right.$$
which can be shortly written as
$$dZ=A\cdot Z\ \ \ \text{where}\ \ \ 
A=\left(\begin{matrix} 
-{\omega_1\over 2}&-\omega_2\\ \omega_0&{\omega_1\over 2}
\end{matrix}\right)\ \ \ \text{and}\ \ \ 
Z=\left(\begin{matrix} z_1\\z_2\end{matrix}\right)$$
The matrix $A$ may be thought as a differential $1$-form on $M$ 
taking values in the Lie algebra $sl(2,\mathbb C)$ satisfying 
integrability condition $dA+A\wedge A=0$.
Then, Darboux Theorem (see \cite{Godbillon}, III, 2.8, iv, p.230)
asserts that there exists, on any simply connected open subset $U\subset M$,
an holomorphic map
$$\Phi:U\to SL(2,\mathbb C)\ \ \ \text{such that}\ \ \ A=\Phi^*\mathcal M$$
where $\mathcal M$ is the Maurer-Cartan $1$-form on $SL(2,\mathbb C)$
(see example \ref{E:SL2}). Moreover, the map $\Phi$ is unique up to composition
by a translation of $SL(2,\mathbb C)$.
\end{remark}

\eject

\begin{example}\label{E:SL2quotient}\rm
Consider the quotient $M:=\Gamma\diagdown SL(2,\mathbb C)$
by a co-compact lattice $\Gamma\subset SL(2,\mathbb C)$.
The left-invariant $1$-forms $(\omega_0,\omega_1,\omega_2)$
defined in example \ref{E:SL2} are well-defined on $M$
and $M$ is parallelizable. Following \cite{HuckleMargu},
there is no non constant meromorphic function on $M$ 
(i.e. the algebraic dimension of $M$ is $a(M)=0$).
Therefore, any foliation $\mathcal F$ on $M$ is defined 
by a global meromorphic $1$-form 
$$\omega=\alpha\omega_0+\beta\omega_1+\gamma\omega_2$$
and the coefficients are actually constants $\alpha,\beta,\gamma\in\mathbb C$.
\end{example}

\begin{cor}
Any foliation $\mathcal F$ on a quotient 
$M:=\Gamma\diagdown SL(2,\mathbb C)$
by a co-compact lattice $\Gamma$ is actually
defined by a left-invariant $1$-form.
In particular, $\mathcal F$ is regular, transversely projective
and minimal: any leaf of $\mathcal F$ is dense in $M$. 
The set of foliations on $M$ is a rational curve.
\end{cor}

%We know, following Lie, that any finite dimensional transitive Lie algebra 
%of germs of vector fields in one complex variable is isomorphic 
%to one of the three following algebras:
%\begin{enumerate}
%\item $\mathbb C\cdot\partial_z$ (euclidean);
%\item $\mathbb C\cdot\partial_z+\mathbb C\cdot z\partial_z$ (affine);
%\item $\mathbb C\cdot\partial_z+\mathbb C\cdot z\partial_z+\mathbb C\cdot z^2\partial_z=sl(2,\mathbb C)$ (projective).
%\end{enumerate}
%They respectively correspond to the translations, affine or Moebius 
%transformation groups of the Riemann sphere.
%Thus there exists three possible transverse structures to 
%a codimension one holomorphic foliation:
%euclidean (and its linear variant), affine or projective structures.

A foliation $\mathcal F$ is {\it transversely euclidean} if there exists 
an atlas of submersions $f_i:U_i \to \mathbb C$ on $M$
defining $\mathcal F$ such that on any $U_i\cap U_j$ we have
\[
f_i = f_j + a_{ij}, \ \ \ a_{ij} \in \mathbb C \, .
\]
Of course, we can {\it glue} the $df_i$ and produce 
a global closed holomorphic $1$-form $\omega_0$ 
inducing $\mathcal F$. In particular $l(\mathcal F)=0$.
By the same way, $\mathcal F$ is {\it transversely linear}
when it can be defined by submersions
$f_i:U_i \to \mathbb C^*$ satisfying the cocycle condition:
\[
f_i = \lambda_{ij} \cdot f_j,  \ \ \  \lambda_{ij} \in \mathbb C^* \, .
\]
Again, we can {\it glue} the  $\frac{df_i}{f_i}$ and 
produce a global closed holomorphic $1$-form 
inducing $\mathcal F$ and we have $l(\mathcal F)=0$.
Via the exponential map, this notion is equivalent to the previous one 
(in the complex setting).

Finally, a foliation $\mathcal F$ is {\it transversely affine} 
when it can be defined by submersions
$f_i:U_i \to \mathbb C$ satisfying the cocycle condition:
\[
f_i = a_{ij} f_j + b_{ij}, \ \ \  a_{ij} \in \mathbb C^*, \,  b_{ij} \in \mathbb C.
\]
Equivalently, an affine structure is locally defined by a pair
of holomorphic $1$-forms $(\omega_0,\omega_1)$ satisfying
$$\left\{\begin{matrix}
d\omega_0 &=& \omega_0 \wedge \omega_1 \, \\
d\omega_1 &=& 0 
\end{matrix} \right. \ \ \ \text{up to modification}\ \ \ 
\left\{\begin{matrix}
\tilde\omega_0 &=& f\cdot\omega_0 \hfill \, \\
\tilde\omega_1 &=& \omega_1-{df\over f}
\end{matrix} \right.$$

%It is not clear whether we can glue the $df_i$: it will depend on the class of $\{ a_{ij} \}$ 
%in ${\rm H}^1(M,\mathbb C^*)$.
%Nevertheless, when $\mathcal F$ is defined by a global $1$-form $\omega_0$, 
%we have that
%\[
%  {\omega_0}\vert_{U_i} = g_i df_i, \ \ \  g_i \in \mathcal O(U_i) \, .
%\]
%Thus, we can write on $U_i \cap U_j$
%\[
%\displaystyle{  {\omega_0}\vert_{U_i\cap U_j} = g_i df_i = a_{ij} g_i df_j = g_j df_j,}
%\]
%which implies that $g_i / g_j = a_{ij}$. In particular, we can define 
%the global closed meromorphic $1$-form $\omega_1$  by
%${\omega_1}\vert_{U_i}= - \frac{dg_i}{g_i}$. Since
%\begin{equation*}
%\displaystyle{\left\lbrace
%\begin{matrix}
%d\omega_0 &=& \omega_0 \wedge \omega_1 \, \\
%d\omega_1 &=& 0 
%\end{matrix} \right. }
%\end{equation*}
%we conclude that $l(\mathcal F) \le 1$.

%Note that, in the euclidean, linear or affine cases, replacing the target of the submersion 
%by the full Riemann sphere $\mathbb CP(1)$ would introduce poles in the $\omega_i$. 

\eject

\subsection{Transversely projective foliations: the singular case \cite{Scardua}}\label{S:projectifSingular}

A singular foliation $\mathcal F$ on a complex manifold $M$ 
will be said {\it transversely projective} if it admits 
a Godbillon-Vey sequence of length $2$, i.e. if there exist
meromorphic $1$-forms $\omega_0$, $\omega_1$ and $\omega_2$ 
on $M$ satisfying $\mathcal F=\mathcal F_{\omega_0}$ and 
$$\left\{\begin{matrix}
d\omega_0 &=& \hfill\omega_0 \wedge \omega_1 \, \\
d\omega_1 &=& 2\omega_0 \wedge \omega_2 \, \\
d\omega_2 &=& \hfill\omega_1 \wedge \omega_2 \,
\end{matrix} \right. $$
The foliation $\mathcal F$ is actually regular and transversely projective 
in the classical sense of \S \ref{S:projectifRegular} on the Zariski
open subset $U=M\setminus((\Omega)_\infty\cup\mathcal Z_0)$ complementary to the set
$(\Omega)_\infty$ of poles for $\omega_0$, $\omega_1$ and $\omega_2$
and the set $\mathcal Z_0$ of zeroes for $\omega_0$ that are not in $(\Omega)_\infty$. 
In fact, $(\omega_0,\omega_1,\omega_2)$ is a regular projective triple 
on $U$. Another triple $(\tilde\omega_0,\tilde\omega_1,\tilde\omega_2)$
defines the same projective structure (on a Zariski open subset) if it is obtained
from the previous one by a combination of
\begin{equation}\label{E:meromorphicTriple}
\left\{\begin{matrix}
\tilde\omega_0 &=& {1\over f}\cdot\omega_0 \hfill \\
\tilde\omega_1 &=& \omega_1 +{df\over f} \\
\tilde\omega_2 &=& f\cdot\omega_2 \hfill
\end{matrix} \right. \ \ \ \text{and}\ \ \ 
\left\{\begin{matrix}
\tilde\omega_0 &=& \omega_0 \hfill \\
\tilde\omega_1 &=& \omega_1  +g\cdot\omega_0\hfill\\
\tilde\omega_2 &=& \omega_2 +g\cdot\omega_1+g^2\cdot\omega_0-dg
\end{matrix} \right. 
\end{equation}
where $f,g$ denote meromorphic functions on $M$.

We note that any pair $(\omega_0,\omega_1)$ satisfying
$d\omega_0=\omega_0\wedge\omega_1$ can be completed 
into a triple subjacent to the projective struture 
in an unique way. It follows that, in the pseudo-parallelizable case, 
a projective transverse structure is always defined by a global
meromorphic triple.

We say that $\mathcal F$ is {\it transversely affine}
if it admits a Godbillon-Vey sequence of length $1$, 
i.e. meromorphic $1$-forms $\omega_0$ and $\omega_1$ 
satisfying
$$\left\{\begin{matrix}
d\omega_0 &=& \omega_0 \wedge \omega_1 \, \\
d\omega_1 &=& 0\hfill 
\end{matrix} \right. $$ 
Another pair $(\tilde\omega_0,\tilde\omega_1)$ will define the same
affine structure if we have
$$\left\{\begin{matrix}
\tilde\omega_0 &=& {1\over f}\cdot\omega_0 \hfill \\
\tilde\omega_1 &=& \omega_1 +{df\over f}
\end{matrix} \right.$$
for a meromorphic function $f$. Finally, we say that $\mathcal F$
is {\it transversely euclidean} (resp. {\it transversely trivial}) 
if it is defined by a closed meromorphic $1$-form $\omega_0$ 
(resp. by an exact $1$-form $\omega_0=df$, $f\in\mathcal M(M)$).

\eject

The foliation $\mathcal H$ defined on $M \times \mathbb CP(1)$
by the integrable $1$-form 
$$\Omega=dz+\omega_0+z\omega_1+z^2\omega_2$$
coincides over $U$ with the suspension of the projective structure,
and will be still called {\it suspension of $\mathcal F$}. In fact,
the vertical hypersurface $(\Omega)_\infty\times\mathbb CP(1)$
is invariant by the foliation $\mathcal H$. Outside of this vertical invariant set,
the foliation $\mathcal H$ is transversal to the vertical $\mathbb CP(1)$-fibration.
Along $\mathcal Z_0$, the foliation $\mathcal H$ is tangent to the zero-section
$M\times\{z=0\}$ and the projective structure ramifies: it is locally defined 
by an holomorphic map $f_i:U_i\mapsto\mathbb CP(1)$ up to composition by an element
of $PGL(2,\mathbb C)$. This ramification set
$\mathcal Z_0$ is invariant for $\mathcal F$ (union of leaves and singular points).
As in the regular case, one can define the monodromy representation
$$\rho:\pi_1(M\setminus(\Omega)_\infty)\to PGL(2,\mathbb C)$$
(ramification points $\mathcal Z_0$ have no monodromy).

In contrast with the regular case, the suspension $\mathcal H$ is well-defined 
only up to a bimeromorphic transformation preserving the generic vertical fibres 
$\{p\}\times \mathbb CP(1)$ and the zero-section $M\times\{z=0\}$
$$\Phi:M \times \mathbb CP(1)\dashrightarrow M \times \mathbb CP(1)\ ;\ 
(p,z)\mapsto(p,f(p) z/(1-g(p)z)),$$
where $f,g\in\mathcal M(M)$ are meromorphic. Note that some irreducible components 
of $(\Omega)_\infty$ may disappear after such a transformation $\Phi$. For instance,
one can show that any irreducible component of $(\Omega)_\infty$ 
which is not $\mathcal F$-invariant may be deleted by a change of triple.
Only the remaining persistent components can generate non trivial local monodromy
for the representation $\rho$. This leads to the following 

\begin{prop}\label{P:singularSimplyConnected}
Let $\mathcal F$ be a (singular) transversely projective (resp. affine) foliation 
on a simply connected manifold $M$. If $(\Omega)_\infty$ has no persistent
component, then $\mathcal F$ admits a meromorphic (resp. holomorphic) first integral.
\end{prop}

\begin{proof}The assumption just means that there exists a covering
$U_i$ of $M$ by Zariski open subset on which the projective structure
can be defined by an holomorphic triple. Therefore, like in 
Remark \ref{R:simplyConnected} the developing map provides a well-defined
meromorphic first integral $f:M\to\mathbb CP(1)$ (possibly with ramifications).
\end{proof}

\begin{cor}\label{C:singularSimplyConnected}Let $\mathcal F$ be a transversely projective 
(resp. affine) foliation on a simply connected manifold $M$. Then 
\begin{itemize}
\item either $\mathcal F$ has a meromorphic (resp. holomorphic) first integral,
\item or $\mathcal F$ admits an invariant hypersurface.
\end{itemize}
\end{cor}

\eject

\begin{remark}\rm
A transversely projective foliation $\mathcal F$ on $M$
with suspension $\mathcal H$ on $M\times\mathbb CP(1)$ is actually 
transversely affine if, and only if, there is a section 
$g:M\to M\times\mathbb CP(1)$ which is invariant by $\mathcal H$.
Indeed, after change of coordinate $\tilde z=z/(1-{z\over g})$
on $M\times\mathbb CP(1)$, we have sent the invariant hypersurface 
$g(M)$ onto $\{z=\infty\}$ which means that $\omega_2=0$.
In the regular case, this is still true after replacing
$M\times\mathbb CP(1)$ by the locally trivial $\mathbb CP(1)$-bundle 
$\pi:N\to M$ (see Proposition \ref{P:Suspension}) and
if we ask moreover that the section $g:M\to N$ has no intersection
with the section $f:M\to N$ providing the projective structure.
\end{remark}

\begin{example}[The Riccati equation over a curve]\label{E:Riccati}\rm
Given meromorphic $1$-forms $\alpha,\beta,\gamma$ on a curve $C$,
the Riccati differential equation
$$dz+\alpha +\beta z+\gamma z^2=0$$
defines a transversely projective foliation $\mathcal H$ on $C\times\mathbb CP(1)$
with meromorphic projective triple
$$\left\{\begin{matrix}
\omega_0 &=&  dz+\alpha +\beta z+\gamma z^2\\
\omega_1 &=&  \hfill \beta+2\gamma z \\
\omega_2 &=&  \hfill \gamma
\end{matrix} \right. $$
The polar set $(\Omega)_\infty$ is the union of the vertical lines 
over the poles of $\alpha,\beta,\gamma$ and the horizontal line
$L_\infty=\{z=\infty\}$. In the chart $w=1/z$, the alternate triple 
$$\left\{\begin{matrix}
\tilde\omega_0 &=&  -dw+\alpha w^2 +\beta w+\gamma \\
\tilde\omega_1 &=&  \hfill -\beta-2\alpha w \\
\tilde\omega_2 &=&  \hfill -\alpha
\end{matrix} \right. $$
(obtained by setting successively $f=1/w^2$ and $g=-2/w$ in (\ref{E:meromorphicTriple}))
shows that $L_\infty$ is not a persistent
pole for the projective structure. When $\gamma=0$, the foliation $\mathcal H$
is transversely affine with poles like above, 
but additionally $L_\infty$ is a persistent zero for the affine structure
(the transverse affine coordinate has a pole along $L_\infty$).

The Riccati foliation above can be thought as the suspension of a singular
projective structure on the curve $C$ (i.e. a dimension $0$ transversely
projective foliation on $C$).
\end{example}

\eject

In the spirit of Theorem \ref{T:CS}, 
one can find in \cite{Scardua} the following

\begin{prop}[Sc\'ardua]Let $\mathcal F$ be a transversely projective foliation
defined by a global meromorphic triple $(\omega_0,\omega_1,\omega_2)$ on $M$.
Assume that the foliation $\mathcal G$ defined 
by $\omega_2$ admits a meromorphic first integral $f\in\mathcal M(M)$.
Then, $\mathcal F$ is the pull-back by a meromorphic map 
$\Phi:M \dashrightarrow C \times \mathbb{C}P(1)$ of the foliation $\mathcal H$ 
defined by a Riccati equation on a curve $C$.
\end{prop}

\begin{proof}One can assume that $\omega_2=df$. Integrability conditions yield
$$  \left\{\begin{matrix}
0=d\omega_2 &=& \hfill \omega_1\wedge\omega_2 \\
\hfill d\omega_1 &=&  2\omega_0\wedge\omega_2 \\
\hfill d\omega_0 &=&  \hfill \omega_0\wedge\omega_1
\end{matrix} \right.\ \ \ \Rightarrow\ \ \  \left\{\begin{matrix}
\omega_1 &=& gdf \\
\omega_0 &=&  {1\over 2}dg+hdf \\
0 &=&  d(h-g^2)\wedge df
\end{matrix} \right.$$
for meromorphic functions $g,h$ on $M$. It follows from Stein Factorization
Theorem that there exists some holomophic map $\phi:M\mapsto C$ onto a curve $C$
through which we can factorize $h-g^2=\tilde h(\phi)$ and $f=\tilde f(\phi)$.
Therefore
$$\omega_0= {1\over 2}dg+\{\tilde h(\phi)+g^2\}\phi^*d\tilde f$$
and $\mathcal F$ is the pull-back via the map $\Phi=(\phi,g)$
of the foliation defined by the Riccati equation
$dz+\tilde hd\tilde f+z^2d\tilde f$.
\end{proof}

\begin{lem}If a foliation $\mathcal F$ admits $2$ distinct 
projective (resp. affine, euclidean) structures, 
then it is actually transversely affine (resp. euclidean, trivial).
\end{lem}

\begin{proof}Assume we have $2$ projective triples $(\omega_0,\omega_1,\omega_2)$
and $(\tilde\omega_0,\tilde\omega_1,\tilde\omega_2)$ that are not related by 
a composition of the admissible changes above: after the admissible change 
setting $\tilde\omega_0=\omega_0$ and $\tilde\omega_1=\omega_1$, we have
$\tilde\omega_2\neq\omega_2$. Therefore, by comparing the second line 
of integrability conditions for both triples, 
we see that $\tilde\omega_2=\omega_2+f\omega_0$ 
for a meromorphic function $f\in\mathcal M(M)$. 
Then, by comparing the third condition, we obtain
$$d(f\omega_0)=\omega_1\wedge(f\omega_0)\ \ \ \text{and thus}\ \ \ 
\omega_0\wedge\omega_1=\omega_0\wedge{df\over 2f}$$
which proves that the pair $(\tilde\omega_0,\tilde\omega_1):=(\omega_0,{df\over 2f})$
is an affine structure for $\mathcal F$. Notice that ${\omega_0\over\sqrt{f}}$
is closed: $\mathcal F$ becomes transversely euclidean on a $2$-fold 
ramified covering of $M$. By the same way, if $(\omega_0,\omega_1)$ and
$(\tilde\omega_0,\tilde\omega_1)$ are $2$ distinct affine structures,
then we may assume $\tilde\omega_0=\omega_0$ and $\tilde\omega_1=\omega_1+f\omega_0$
with $d\omega_1=d(f\omega_0)=0$ and conclude that $\mathcal F$ is actually
defined by the closed meromorphic $1$-form $f\omega_0$. Finally,
if $\omega_0$ and $f\omega_0$ are $2$ closed meromorphic $1$-forms
defining $\mathcal F$, then $f$ is a meromorphic first integral for $\mathcal F$.
\end{proof}

\eject

The present singular notion of transversely projective foliation
is clearly stable under bimeromorphic transformations. Moreover,
the main result of \cite{Casale} permits to derive

\begin{thm}\label{T:Casale}
Let $\phi:\tilde M\dashrightarrow M$ be a dominant meromorphic map 
between compact manifolds and let $\mathcal F$ 
be a foliation on $M$ admitting a Godbillon-Vey sequence.
Then, $\tilde{\mathcal F}=\phi^*\mathcal F$ is transversely projective 
(resp. affine) if, and only if, so is $\mathcal F$.
\end{thm}

The analogous result for transversely euclidean foliations is false:
one can find in \cite{LinsNeto} an example of a transversely affine 
foliation which becomes transversely euclidean on a finite covering
(a linear foliation on a torus). The assumption dominant is necessary
since there are examples of non transversely projective foliations
which become transversely affine in restriction to certain non tangent 
hypersurface (see section \ref{S:Alcides}). 

\begin{proof}
%Of course, any (singular) projective, affine or euclidean transverse 
%structure for $\mathcal F$ lift as a similar structure for $\tilde{\mathcal F}$
%without assumption: the pull-back by any non constant meromorphic map 
%$\tilde M\dashrightarrow M$ of a lenght $k$ Godbillon-Vey sequence 
%is also a length $k$ Godbillon-Vey sequence. 
Since a Godbillon-Vey sequence can be pulled-back by any non constant 
meromorphic map, we just have to prove that projective (resp. affine) 
structure can be pushed-down under the assumptions above.
In the case $\phi$ is a finite ramified covering, then the statement 
is equivalent to Theorem 1.6 (resp. 1.4) in \cite{Casale}.

In the case $\phi$ is holomorphic with connected generic fibre, 
then choose meromorphic $1$-forms $\omega_0$
defining $\mathcal F$ and $\omega_1$ satisfying 
$d\omega_0=\omega_0\wedge\omega_1$ on $M$ and consider their pull-back
$\tilde\omega_0$ and $\tilde\omega_1$ on $\tilde M$.
Then, there is a unique meromorphic $1$-form $\tilde\omega_2$
completing the previous ones into a projective triple
compatible with the structure of $\tilde{\mathcal F}$.
On the other hand, reasonning as in Lemma \ref{L:tangente}
at the neighborhood $\tilde U=\phi^{-1}(U)$ of a generic fibre $\phi^{-1}(p)$,
we see that the foliation $\tilde{\mathcal F}$ is defined by a submersion 
$\tilde f:\tilde U\mapsto\mathbb CP(1)$ defining the projective structure
and can be pushed-down into a submersion $f:U\mapsto\mathbb CP(1)$.
This latter one defines a projective structure transverse 
to $\mathcal F$ on $U$. There exists a unique meromorphic $1$-form 
$\omega_2$ on $U$ completing $\omega_0$ and $\omega_1$ into
a compatible projective triple. By construction, $\tilde\omega_2$
must coincide with $\phi^*\omega_2$ on $\tilde U$. Therefore,
$\tilde\omega_2$ is tangent to the fibration given by $\phi$ on $\tilde U$,
and thus everywhere on $\tilde M$.
By connexity of the fibres, $\tilde\omega_2$ is actually the pull-back 
of a global meromorphic $1$-form $\omega_2$ on $M$ (which extends 
the one previously defined on $U$).

Finally, by Stein Factorization Theorem, the statement reduces to
the two cases above.
\end{proof}

\eject
%%%%%%%%%%%%%%%%%%%%%%%%%%%%%%%%%%%%%%%%%%%%%%%%%%%%%%%%%%%%%%%%%%%%%%%%%%%%%%%%%%%%%%%%%%%%%%%%%%%%%%%%%%%%%%%%%%%%%%%%%%%%%%%%

\section{Proof of Theorem \ref{T:Main}}\label{S:main}

Let $\mathcal F$ be a foliation on a compact manifold $M$ 
admitting a Godbillon-Vey sequence $(\omega_0,\omega_1,\ldots)$
and consider the maximal non trivial form $\Theta=\omega_0\wedge\ldots\wedge\omega_{n-1}$:
$\Theta\wedge\omega_{n}=0$. Like in the introduction, 
we denote by $K$ the field of meromorphic first integrals for $\Theta$
and consider the reduction map $\pi:M\dashrightarrow N$ associated to this field.
The fibration $\mathcal G$ induced by $\pi$ contains the foliation 
$\mathcal F_\Theta$ as a sub-foliation, and may be of larger dimension 
as soon as there are few meromorphic functions on $M$. In particular,
there is no reason why $\mathcal G$ is a sub-foliation of $\mathcal F$. 
Anyway, when $\mathcal G\subset\mathcal F$, then we are in the first alternative of
Theorem \ref{T:Main}: {\it $\mathcal F$ is the pull-back 
of an algebraic foliation $\underline{\mathcal F}$ on $\text{red}(M,\Omega)$}.
Indeed, after modification of $M$, one can assume that 
the reduction map is holomorphic with connected fibres. 
The claim above immediately follows from:

\begin{lem}\label{L:tangente}
Let $\mathcal F$ be a foliation on a complex manifold $M$.
Let  $\pi: M \to N$ be a surjective holomorphic map whose fibers
are connected and tangent to $\mathcal F$, 
that is, contained in the leaves of $\mathcal F$. 
Then, $\mathcal F$ is the pull-back by $\pi$ of a foliation
$\tilde{\mathcal F}$ on $N$.
\end{lem}

\begin{proof}
In a small connected neighborhood $U\subset M$ of a generic point $p \in M$, 
the foliation $\mathcal F$
is regular, defined by a local submersion $f:U\to\mathbb C$. Since $f$ is contant
along the fibers of $\pi$ in $U$, we can factorize $f=\tilde f\circ\pi$ for
an holomorphic function $\tilde f:\pi(U)\to\mathbb C$. In particular, the function
$\tilde f$ defines a codimension one singular foliation $\tilde{\mathcal F}$
on the open set $\pi(U)$. Of course, $\tilde{\mathcal F}$ does not depend
on the choice of $f$. Moreover, since $f=\tilde f\circ\pi$, 
the function $f$ extends to the whole tube $T:=\pi^{-1}(\pi(U))$.
By connectivity of $U$ and the fibers of $\pi$, the tube $T$ is connected 
and the foliation $\mathcal F$ is actually defined by $f$ on the whole of $T$,
coinciding with $\pi^*(\tilde{\mathcal F})$ on $T$.
In this way, we can define a foliation $\tilde{\mathcal F}$ on
$N\setminus S$, where $S=\{p\in N\ ;\ \pi^{-1}(p)\subset\text{Sing}(\mathcal F)\}$
such that $\mathcal F=\pi^*(\tilde{\mathcal F})$. We note that $S$
has codimension $\ge2$ in $N$; therefore, $\tilde{\mathcal F}$ extends on $N$
by Levy's Extension Theorem.
\end{proof}
 
We now assume that the fibration $\mathcal G\not\subset\mathcal F$.
We note
that when $n=\dim(M)$, then $M$ is actually pseudo-parallelizable
and the field $K$ coincides with the field $\mathcal M(M)$ of 
meromorphic functions on $M$. 

\eject

We introduce the sheaf $\mathcal B$ of basic meromorphic vector fields for $\mathcal F_\Theta$; 
a section $X$ of $\mathcal B(U)$ over $U\subset M$ is characterized by the following property:
\begin{equation}\label{E:DefinitionBasic}
L_X\Omega=d(i_X\Theta)=f\cdot\Theta,\ \ \ \text{for some}\ f\in\mathcal M(U)
\end{equation}
(here, we use that $L_X=d\circ i_X+i_X\circ d$ and the fact that $d\Theta=0$).
We remark that $\mathcal B$ is a sheaf of Lie algebras over $\mathbb C$. 
The subsheaf of vector fields tangent to $\mathcal F_\Theta$ 
\begin{equation}\label{E:DefinitionI}
\mathcal I(U):=\{X\in\mathcal B(U)\ ;\ i_X\Theta=0\}.
\end{equation} 
form a Lie-ideal of $\mathcal B$: if $X\in\mathcal B(U)$ and $Y\in\mathcal I(U)$,
then  $[X,Y]\in\mathcal I(U)$ as can be seen from a local flow-box for $\mathcal F_\Theta$.
The quotient $\mathcal T=\mathcal B/\mathcal I$ is a sheaf of Lie algebras over $\mathbb C$,
whose sections are the {\bf transversal relative vector fields} to $\mathcal F_\Theta$.
Although the sheaf $\mathcal B$ may have no global meromorphic section, the relative sheaf 
$\mathcal T$ has many, as shown by the:

\begin{lem}\label{L:BasisTransversal}
Let $T:=\mathcal T(M)$ be the Lie algebra over $\mathbb C$ of global
transversal relative vector fields. Then $T$ is a $n$-dimensional vector space over $K$
and admits a canonical basis $(\mathbf X_0,\ldots,\mathbf X_{n-1})$ satisfying 
$$\omega_k(\mathbf X_l)=\delta^k_l\ \ \ \ \text{for}\ \ \ k,l=0,\ldots,n-1.$$
\end{lem}

We note that $\omega_k(\mathbf X_l)$ is well defined since $\mathbf X_l$ 
is locally defined
modulo an element $Y\in\mathcal I$ and $\omega_k(Y)=0$. 
More generally, we will use the fact that an element $\mathbf X\in T$
acts as a derivation on $K$: $\mathbf X\cdot f\in K$ for all $f\in K$.
Indeed, for any local representative $X$ of $\mathbf X$, 
$X$ is a basic vector field and $X\cdot f$ is a local first integral 
for $\mathcal F_\Theta$; since $Y\cdot f=0$ for any $Y\in\mathcal I$,
we can set unambiguously $\mathbf X\cdot f:=X\cdot f$ which is now a global
first integral, thus belonging to $K$. Similarly, $L_{\mathbf X}\omega_k$
is a well-defined global meromorphic $1$-form on $M$ for any $\mathbf X\in T$
and $k=0,\ldots,n-1$.

\eject

Before proving the Lemma,
we note that, by maximality property of $n$, one can write 
\begin{equation}\label{E:DefinitionomegaN+1}
\omega_n=a_1\omega_1+\cdots+a_{n-1}\omega_{n-1},\ \ \ a_k\in\mathcal M(M)
\end{equation}
(here, (\ref{E:changeOmegak}) allow us to set $a_0=0$). 
In fact, coefficients $a_k$ actually belong to $K$.
Indeed, after combining (\ref{E:condition}) and (\ref{E:DefinitionomegaN+1}), we get
\begin{equation}\label{E:DerivationomegaN}
d\omega_{n-1}=\omega_0\wedge\sum_{k=1}^{n-1} a_k\omega_k
+\sum_{k=1}^{n-1}\left(\begin{matrix}k\\n-1\end{matrix}\right)\omega_k\wedge\omega_{n-k}.
\end{equation}
After differentiation and multiplication by the $n-1$-form
\begin{equation}\label{E:DefinitionThetak}
\widehat{\Theta}_k
=\omega_0\wedge\cdots\wedge\omega_{k-1}\wedge\omega_{k+1}\wedge\cdots\wedge\omega_{n-1},\
\ \ k=1,\ldots,n-1,
\end{equation}
we obtain
$$\Theta\wedge da_k=0$$
and, as a consequence, that $a_k$ is a first integral for $\Theta$.

\begin{proof}From a local flow-box for $\mathcal F_\Theta$, one easily see
that $T$ is a $K$-vector space: if $\mathbf X\in T$ and $f\in K$, then $f\cdot\mathbf X\in T$.
Now, consider local vector fields $(X_0,\ldots,X_{n-1})$ dual to $(\omega_0,\ldots,\omega_{n-1})$
like in the statement: they are well-defined modulo sections of $\mathcal I$
and define global relative vector fields $(\mathbf X_0,\ldots,\mathbf X_{n-1})$; 
we have to prove that they are transverse
relative vector fields, i.e. that $L_{\mathbf X_k}\Theta=f\Theta$ for some $f\in\mathcal M(M)$
(actually, $f\in K$ since $f\Theta$ will be closed as well).
We have $L_{\mathbf X_k}\Theta=d(\widehat\Theta_k)$ where $\widehat\Theta_k$
is defined by (\ref{E:DefinitionThetak}).
From Godbillon-Vey relations (\ref{E:condition}), one easily deduce that
all $d(\widehat\Theta_k)$ are zero,
except $d(\widehat\Theta_1)=c\cdot\Theta$ for some constant $c\in\mathbb C$
and $d(\widehat\Theta_0)=\pm\widehat\Theta_{n-1}\wedge\omega_n=\pm a_{n-1}\Theta$.

Now, given an element $\mathbf X\in T$, one can
write $\mathbf X=\lambda_0\mathbf X_0+\cdots+\lambda_{n-1}\mathbf X_{n-1}$ 
modulo $\mathcal I$
where $\lambda_k=\omega_k(\mathbf X)$ are global meromorphic functions; 
as can be seen for a local flow-box, all $\lambda_i$ must be first integrals 
for $\mathcal F_\Theta$.
\end{proof}

The proof of Theorem \ref{T:Main} is similar to  
that of \cite{Crocodile} after substituting global sections
of $\mathcal T$ to global meromorphic vector fields.
We consider the Lie sub-algebra 
$$\mathcal L:=\{\mathbf X\in T\ ;\ \mathbf X\cdot K=0\}$$
of those relative vector fields that are tangent to the fibration $\mathcal G$
given by global first integrals of $\Theta$. We note that $\mathcal L$
is now a Lie algebra over $K$. Assuming that $\mathcal G\not\subset\mathcal F$
(otherwise, we have already concluded the proof by (\ref{L:tangente}),
we consider the Lie sub-algebra (over $K$) defined by
$$\mathcal L_0:=\{\mathbf X\in\mathcal L\ ;\ \omega_0(\mathbf X)=0\}.$$
Clearly, $\mathcal L_0$ is a codimension $\le 1$ sub-algebra of $\mathcal L$
over $K$; we now prove that indeed $\mathcal L/\mathcal L_0$ is not trivial:

\begin{lem}\label{L:ExistenceX}If $\mathcal G\not\subset\mathcal F$, 
there exists $\mathbf X\in T$ such that $\omega_0(\mathbf X)=1$ 
and $\mathbf X\cdot K=0$,
i.e. $\mathbf X(f)=0$ for any $f\in K$.
\end{lem}

\begin{proof} If $K=\mathbb C$, then the Lemma is trivial.
If not, suppose that $f_1,\ldots,f_{N}$ are elements of $K$
such that 
$$df_1\wedge\cdots\wedge f_{N}\not\equiv 0$$
with $N$ maximal: we have by asumption $N<n$ and
\begin{equation}\label{E:TechnicalLemma}
\omega_0\wedge df_1\wedge\cdots\wedge df_{N}\not\equiv 0.
\end{equation}
Remark now that if $\mathbf X\in T$ satisfies $\mathbf X(f_k)=0$ for $k=1,\ldots,N$,
then for each $f\in K$, the meromorphic function $\mathbf X(f)$
is actually zero. Let us write
$$\mathbf X=\alpha_0\mathbf X_0+\alpha_1 \mathbf X_1+\cdots
+\alpha_{n-1} \mathbf X_{n-1}$$
with $\alpha_k\in K$; since $\omega_0(\mathbf X)=\alpha_0$, we already set
$\alpha_0:=1$. We now have to solve the $N\times(n-1)$-linear system
$$\left\{\begin{matrix}
df_1(\mathbf X_0)+df_1(\mathbf X_1)\cdot\alpha_1+\cdots
+df_1(\mathbf X_{n-1})\cdot\alpha_{n-1}&=&0\\
\vdots&=&0\\
df_N(\mathbf X_0)+df_N(\mathbf X_1)\cdot\alpha_1+\cdots
+df_N(\mathbf X_{n-1})\cdot\alpha_{n-1}&=&0
\end{matrix}\right.$$
From (\ref{E:TechnicalLemma}), the corresponding matrix $(df_k(\mathbf X_l))_{k,l}$
has maximal rank and one can solve the system above. 
If $N<n-1$, there are obviously many solutions.
\end{proof}

The proof of Lemma 3.2 in \cite{Crocodile} may be transposed 
to our relative setting:

\begin{lem}\label{L:TechnicalCrocodile}
If the relative vector field $\mathbf X\in T$ satisfies 
$\omega_0(\mathbf X)=1$ and $\mathbf X\cdot f=0$ for any $f\in K$
like in Lemma \ref{L:ExistenceX}, then 
$$(L^{(k)}_{\mathbf X}\omega_0)(\mathbf Y)
=(-1)^k \omega_0(L^{(k)}_{\mathbf X}(\mathbf Y))$$
for any $\mathbf Y\in T$; here, we denote by $L_X(Y)$ the Lie bracket $[X,Y]$.
\end{lem}

\eject

Now, we can assume that the $\omega_k$ 
are given by $\omega_k=L^{(k)}_{\tilde{\mathbf X}_0}\omega_0$: 
this modification does not affect
neither the foliation $\mathcal F_\Theta$, nor the field $K$.
We keep on notations $T=K<\mathbf X_0,\ldots,\mathbf X_{n-1}>$, 
$\mathcal L$, $\mathcal L_0$, etc...
We are going to prove that, after conveniently choosing the generator
$\mathbf X$ for $\mathcal L/\mathcal L_0$ given by Lemma \ref{L:ExistenceX}, 
then we have $\omega_3=L^{(3)}_{\mathbf X}\omega_0\equiv0$
and $\mathcal F$ is transversely projective, thus concluding the proof of 
Theorem \ref{T:Main}.

Given a Lie algebra $\mathcal L$ over a field $K$ of characteristic $0$
and a codimension $1$ subalgebra $\mathcal L_0$ like above, 
a result due to J. Tits (see \cite{Tits}, or \cite{Crocodile0}, p.31-33) 
asserts that there exists an ideal $\mathcal J\subset\mathcal L_0$ 
having codimension $\le 3$ in $\mathcal L$ 
and the quotient $\mathcal L/\mathcal J$ is of one of the following three types:
\begin{enumerate}
\item $\mathcal L/\mathcal J\simeq K\cdot \mathbf X$ 
and $\mathcal J=\mathcal L_0$,
\item $\mathcal L/\mathcal J\simeq K\cdot \mathbf X + K\cdot \mathbf Y$ 
with $[\mathbf  X,\mathbf  Y]=\mathbf  X$ 
and $\mathcal L_0/\mathcal J=K\cdot \mathbf Y$,
\item $\mathcal L/\mathcal J\simeq K\cdot\mathbf X + K\cdot\mathbf Y 
+ K\cdot\mathbf Z$ with $\text{sl}(2)$ relations
$$[\mathbf X,\mathbf Y]=\mathbf X,\ \ \ [\mathbf X,\mathbf Z]=2\mathbf Y\ \ \ \text{and}\ \ \ 
[\mathbf Y,\mathbf Z]=\mathbf Z.$$
\end{enumerate}
In each case, $\mathbf X$ is one of the vector fields produced
by Lemma (\ref{L:ExistenceX}).

In order to prove that $\omega_3=L^{(3)}_{\mathbf X}\omega_0\equiv0$,
we just have to verify that $\omega_3(\mathbf V)=0$ for any relative vector field
$\mathbf V\in T$. Indeed, any local meromorphic vector $V$ field decomposes as 
$$V=f_0\cdot\mathbf X_0+\cdots+f_{n-1}\cdot\mathbf X_{n-1}+V'$$ 
where $f_k$ are local meromorphic functions, 
$\mathbf X_k$ are local representatives for the basis given 
by Lemma \ref{L:BasisTransversal} and $V'$ is a vector field tangent 
to $\Theta$ (and in particular to $\omega_3$); we thus have 
$$\omega_3(V)=f_0\cdot\omega_3(\mathbf X_0)+\cdots+f_{n-1}\cdot\omega_3(\mathbf X_{n-1}).$$
By Lemma \ref{L:TechnicalCrocodile}, we just have to prove that 
$\omega_0(L^{(3)}_{\mathbf X}\mathbf V)=0$ for any $\mathbf V\in T$.
In fact, it is enough to consider $\mathbf V\in T_0$ since $L_{\mathbf X}\mathbf X=0$.
Since $T$ acts by derivation on $K$ and $\mathcal L$ is the kernel,
Observe that $\mathcal L$ is an ideal of $T$: since the elements of $T$ 
act as derivation on $K$, they can be considered as basic vector fields
with regards to the fibration $\mathcal G$ while $\mathcal L$ is the 
sub-algebra of tangent vector fields. In particular, for any $\mathbf V\in T_0$,
we have $L_{\mathbf X}\mathbf V=[\mathbf X,\mathbf V]\in \mathcal L$.

\eject

We now discuss on the three cases given by Tits' Result. 

{\bf First Case: $\mathcal L_0$ is an ideal of $\mathcal L$.}

We have $[\mathbf X,\mathbf V]=f\cdot\mathbf  X$ modulo $\mathcal L_0=\mathcal J$ 
for some $f\in K$. Therefore, $[\mathbf X,[\mathbf X,\mathbf V]]=0$ 
modulo $\mathcal J$ and $\omega_0(L^{(2)}_{\mathbf X}\mathbf V)=0$:
in this case, the foliation $\mathcal F$ is transversely affine.

{\bf Second Case: $\mathcal L_0/\mathcal J$ is generated by $\mathbf Y$ 
with $[\mathbf X,\mathbf Y]=\mathbf X$ modulo $\mathcal J$.}

We have $[\mathbf X,\mathbf V]=f\cdot\mathbf  X+g\cdot\mathbf  Y$ mod $\mathcal J$ 
and $[\mathbf Y,\mathbf V]=h\cdot Y$ mod $\mathcal J$ 
for coefficients $f,g,h\in K$ (here, we use the fact that 
both $\mathbf Y$ and $\mathbf V$ are tangent to $\mathcal F$, whence their Lie bracket).
Applying Jacobi identity to $\mathbf X$, $\mathbf Y$ and $\mathbf V$ yields:
$$[\mathbf X,[\mathbf Y,\mathbf V]]+[\mathbf V,[\mathbf X,\mathbf Y]]
+[\mathbf Y,[\mathbf V,\mathbf X]]=h\cdot\mathbf  X-g\cdot\mathbf  Y=0$$
and we have $h=g=0$. In particular, $[\mathbf X,\mathbf V]=f\cdot\mathbf  X$ 
and $[\mathbf X,[\mathbf X,\mathbf V]]=0$.
We conclude as before that $\mathcal F$ is transversely affine.

{\bf Third Case: $\mathcal L_0/\mathcal J$ is generated by $\mathbf Y,\mathbf Z$ 
with $[\mathbf X,\mathbf Y]=\mathbf X$, $[\mathbf X,\mathbf Z]=2\mathbf Y$ and $[\mathbf Y,\mathbf Z]=\mathbf Z$ modulo $\mathcal J$.}

We have:
$$\left\{\begin{matrix}
[\mathbf X,\mathbf V]&=&f\cdot \mathbf X+g\cdot \mathbf Y+h\cdot \mathbf Z\\
[\mathbf Y,\mathbf V]&=&\hfill i\cdot \mathbf Y+j\cdot \mathbf Z\\
[\mathbf Z,\mathbf V]&=&\hfill k\cdot \mathbf Y+l\cdot \mathbf Z\\
\end{matrix}\right.\ \ \ \text{mod}\ \ \mathcal J$$
for some coefficients $f,g,h,i,j,k,l\in K$. Jacobi identity yields:
$$\begin{matrix}
[\mathbf X,[\mathbf Y,\mathbf V]]+[\mathbf V,[\mathbf X,\mathbf Y]]+[\mathbf Y,[\mathbf V,\mathbf X]]=\hspace{3cm}\\ \\
\hfill =i\cdot \mathbf X+(2j-g)\cdot \mathbf Y-2h\cdot \mathbf Z=0\\ \\
[\mathbf X,[\mathbf Z,\mathbf V]]+[\mathbf V,[\mathbf X,\mathbf Z]]+[\mathbf Z,[\mathbf V,\mathbf X]]=\hfill\\ \\
\hfill =k\cdot \mathbf X+2(f+l-i)\cdot \mathbf Y+(g-2j)\cdot \mathbf Z=0\\ \\
$$[\mathbf Y,[\mathbf Z,\mathbf V]]+[\mathbf V,[\mathbf Y,\mathbf Z]]+[\mathbf Z,[\mathbf V,\mathbf Y]]=\hfill
-k\cdot \mathbf Y+i\cdot \mathbf Z=0\end{matrix}$$
modulo $\mathcal J$ and thus $h=i=k=0$, $l=-f$ and $g=2j$. In particular, 
$[\mathbf X,[\mathbf X,[\mathbf X,\mathbf V]]]=0$
and $\mathcal F$ is transversely projective,
thus proving the Theorem \ref{T:Main}. 

\eject

\section{Proof of Theorem \ref{T:CS}}\label{S:Alcides}

In fact, we prove the more precise

\begin{thm}\label{T:CSprecise}
Let $\mathcal F$ be a foliation admitting a finite Godbillon-Vey sequence 
$(\omega_0,\omega_1,\ldots,\omega_N)$ of length $N+1\ge4$.
Then
\begin{itemize}
\item either $\mathcal F$ is the pull-back by a meromorphic map 
$\Phi:M \dashrightarrow C \times \mathbb{C}P(1)$ of the foliation 
$\underline{\mathcal F}$ defined by 
$$dz+\underline{\omega}_0+\underline{\omega}_1z+\cdots+\underline{\omega}_Nz^n$$
where $\underline{\omega}_k$ are meromorphic $1$-forms on the curve $C$,
\item or $\mathcal F$ is transversely affine.
\end{itemize}
\end{thm}

In particular, we see that a purely transversely projective foliation 
cannot admits other finite Godbillon-Vey sequences than the projective 
triples. 

\begin{proof}
Following Lemma \ref{L:tangGVfini}, we have
$$\Omega=dz+\omega_0+z\omega_1+\left(\sum_{k=2}^N f_k\cdot z^k\right)\omega_N$$
for meromorphic functions $f_k\in\mathcal M(M)$, $f_N\equiv 1$ and
$\omega_N\not=0$.
If $f_{N-1}=0$, then integrability conditions imply that $d\omega_1=0$
(see Lemma \ref{L:tangGVfini}) and $\mathcal F$ is transversely affine. Otherwise,
after a change of Godbillon-Vey sequence of the form (\ref{E:changeOmega0})
(see Section \ref{S:GodbillonVey}), we may assume moreover $f_{N-1}=N$.
Now, the change of coordinate $\tilde z=z+1$ on $\Omega$
$$\begin{matrix}
\Omega&=&d(\tilde z-1)+\omega_0+(\tilde z-1)\omega_1+\cdots+(\tilde z-1)^N\omega_N\\
&=&d\tilde z+\tilde\omega_0+\tilde z\tilde\omega_1+\cdots+\tilde z^N\tilde\omega_N\hfill
\end{matrix}$$
provides a new sequence $(\tilde\omega_0,\tilde\omega_1,\ldots,\tilde\omega_N)$ 
of length $N+1$ satisfying integrability conditions (\ref{E:condition}) (see Introduction).
We take care that this is not a new Godbillon-Vey sequence for $\mathcal F$
(but for $\mathcal F_{\tilde\omega_0}$, whenever $\tilde\omega_0\not=0$).
In fact, we have 
$$\omega_0=\tilde\omega_0+\tilde\omega_1+\tilde\omega_2+\cdots+\tilde\omega_N.$$
We also note that $\tilde\omega_N=\omega_N$ and 
$\tilde\omega_{N-1}=0$.
Following Lemma \ref{L:tangGVfini}, 
there exist meromorphic functions $g_k$ satisfying
\begin{equation}\label{E:CS1}
\tilde\omega_k=g_k\cdot\omega_N\ \ \ \ 
\text{for}\ \ \ k=0,2,\ldots,N-2
\end{equation}
and integrability conditions now write
\begin{equation}\label{E:CS2}
d\tilde\omega_k=(k-1)\tilde\omega_1\wedge\tilde\omega_k\ \ \ \ 
\text{for}\ \ \ k=0,2,\ldots,N-2
\end{equation}
and
\begin{equation}\label{E:CS3}
d\omega_N=(N-1)\tilde\omega_1\wedge\omega_N,\ \ \ d\tilde\omega_1=0.
\end{equation}

\eject

In particular, we see that $\omega_N$ is transversely affine and that
\begin{equation}\label{E:CS4}
\omega_0=\tilde\omega_1+(g_0+g_2+\cdots+g_{N-2})\omega_N.
\end{equation}
Following Lemma \ref{L:technical1} below, there is a non constant meromorphic
function $g\in\mathcal M(M)$ such that $dg\wedge\omega_N=0$. It follows from
Stein's Factorization Theorem that there exist:
\begin{itemize}
\item a meromorphic map $\phi:M\dashrightarrow C$ onto a smooth, compact, 
complex and connected curve $C$,
\item a meromorphic function $\underline g:C\dashrightarrow\overline{\mathbb C}$,
\end{itemize}
such that $g=\underline g\circ\phi$ and the generic fibers $\phi^{-1}(c)$
are irreducible hypersurfaces of $M$. Let $\underline\omega$ be a non zero
meromorphic $1$-form on $C$. The $1$-form $\omega:=\phi^*\underline\omega$
on $M$ is closed, non zero and $df\wedge\omega=0$. Therefore, we can write
$\tilde\omega_N=h_N\cdot\omega$ for a meromorphic function $h$ and setting
$h_k=h_N\cdot g_k$, we get
\begin{equation}\label{E:CS5}
\tilde\omega_k=h_k\omega\ \ \ \ \text{for}\ \ \ k=0,2,\ldots,N-2,N.
\end{equation}
From equations (\ref{E:CS2}), we deduce that
\begin{equation}\label{E:CS6}
\left\{\begin{matrix}\text{either}\ h_k=0,\hfill\\
\text{or}\ 
\left(\tilde\omega_1-{1\over k-1}{dh_k\over h_k}\right)\wedge\omega=0
\end{matrix}\right.
\ \ \ \ 
\text{for}\ \ \ k=0,2,\ldots,N-2,N.
\end{equation}
Thus, for any $k,l=0,2,\ldots,N-2,N$ such that $h_k,h_l\neq0$, we have
\begin{equation}\label{E:CS7}
\left({1\over k-1}{dh_k\over h_k}-{1\over l-1}{dh_l\over h_l}\right)\wedge\omega=0
\end{equation}
and ${h_k^{(l-1)}\over h_l^{(k-1)}}$ is a first integral for $\omega$.
Let $r=gcd\{k-1\ ;\ h_k\neq0\}$: we have
$\sum_{h_k\neq0}n_k(k-1)=r$ for integers $n_k$. Set 
$$h:=\prod_{h_k\neq0}h_k^{n_k}.$$
Therefore, summing equations (\ref{E:CS6}) over $l$, we get
$$
0=\sum_{h_l\neq0}n_l\left({l-1\over r}{dh_k\over h_k}-{k-1\over r}{dh_l\over h_l}\right)\wedge\omega
=\left({dh_k\over h_k}-{k-1\over r}{dh\over h}\right)\wedge\omega.
$$

\eject

Thus, ${h_k\over h^{k-1\over r}}$ is a first integral for $\omega$
and we can write
\begin{equation}\label{E:CS8}
\left\{\begin{matrix}\text{either}\ h_k=0,\hfill\\
\text{or}\ 
h_k=\underline h_k\circ\phi\cdot h^{k-1\over r}
\end{matrix}\right.
\ \ \ \ 
\text{for}\ \ \ k=0,2,\ldots,N-2,N
\end{equation}
and meromorphic functions $\underline h_k:C\dashrightarrow\overline{\mathbb C}$.
From equation (\ref{E:CS4}), we deduce
\begin{equation}\label{E:CS9}
\omega_0=\tilde\omega_1+
\left(\sum_{k=0,2,\ldots,N}\underline h_k\circ\phi\cdot h^{k-1\over r}\right)\omega
\end{equation}
(setting $\underline h_k=0$ whenever $h_k=0$). On the other hand,
from (\ref{E:CS6}) and (\ref{E:CS8}) we get
\begin{equation}\label{E:CS10}
\tilde\omega_1\wedge\omega={1\over k-1}{dh_k\over h_k}\wedge\omega=
{1\over r}{dh\over h}\wedge\omega
\end{equation}
and $\tilde\omega_1={1\over r}{dh\over h}+f\omega$ for a meromorphic function $f$.
Since $\tilde\omega_1$ and $\omega$ are closed, we get after derivating
equation (\ref{E:CS10}) that $df\wedge\omega=0$, i.e. we can write 
$f=\underline f\circ \phi$ for a meromorphic function 
$\underline f:C\dashrightarrow\overline{\mathbb C}$.
Finally, we obtain
$$
\omega_0={1\over r}{dh\over h}+\underline f\circ\phi\cdot\omega
+\left(\sum_{k=0,2,\ldots,N}\underline h_k\circ\phi\cdot h^{k-1\over r}\right)\omega
$$
and, setting $\Phi=(\phi,h):M\dashrightarrow C\times\mathbb CP(1)$
$$
rh\omega_0=\Phi^*\left(dz+r(\underline f\cdot z+
\sum_{k=0,2,\ldots,N}\underline h_k z^{{k-1\over r}+1})\omega\right).
$$ 
\end{proof}

\begin{lem}\label{L:technical1}
Let $\mathcal F$ be a foliation admitting a finite Godbillon-Vey sequence
$(\omega_0,\omega_1,\ldots,\omega_N)$ of length $N+1\ge4$. 
Then
\begin{itemize}
\item either $\omega_N=fdg$ for meromorphic functions $f,g\in\mathcal M(M)$,
\item or $\mathcal F$ is transversely affine.
\end{itemize}
\end{lem}

\begin{proof}We start as in proof of Theorem \ref{T:CSprecise}, keeping
the same notations. Substituting (\ref{E:CS1}) into integrability conditions 
(\ref{E:CS2}) yield
$$(dg_k+(N-k)g_k\tilde\omega_1)\wedge\tilde\omega_N=0\ \ \ \ \text{for}\ \ \ \ k=0,2,\ldots,N-2.
$$

\eject

{\bf If} there exist two distinct integers $k,l\in\{0,2,\ldots,N-2\}$ such that
$g_k,g_l\neq0$, then we can deduce that
$$\left((N-k){dg_l\over g_l}-(N-l){dg_k\over g_k}\right)\wedge\tilde\omega_N=0\ ;$$
{\bf if moreover} the left factor is not zero, then we can conclude that
$$dg\wedge\tilde\omega_N=0\ \ \ \text{with}\ \ \ 
g:={g_l^{(N-k)}\over g_k^{(N-l)}}\ \ \ \text{non constant}$$
i.e. $\omega_N=fdg$ for some meromorphic function $f$. 
{\bf Otherwise}, the discussion splits into many cases.

\noindent{\bf Case 1.} Assume that $g_k=0$ for all $k\in\{0,2,\ldots,N-2\}$.
Then 
$$\omega_0=\sum_{k=0}^N\tilde\omega_k=\tilde\omega_1+\tilde\omega_N$$
and, since $d\tilde\omega_1=0$, we have
$$d\omega_0=d\tilde\omega_N=(N-1)\tilde\omega_1\wedge\tilde\omega_N
=(N-1)\tilde\omega_1\wedge\omega_0$$
and $\mathcal F$ is transversely affine.

%\noindent{\bf Case 2.} Assume now that there is exactly one integer
%$k\in\{0,2,\ldots,N-2\}$ such that $g_k\neq0$. Denote
%$$\beta=\tilde\omega_1+{1\over N-k}{dg_k\over g_k}.$$
%\noindent{\bf Subcase 2.1:} $\beta=0$. Since 
%$\omega_0=\tilde\omega_1+(g_k+1)\tilde\omega_N$,
%we get that either $g_k+1=0$ and $\omega_0=\tilde\omega_1$
%is closed, or $g_k+1\neq0$ and we have
%$$d{\omega_0\over 1+g_k}=d\tilde\omega_N=(N-1)\tilde\omega_1\wedge\tilde\omega_N
%=(N-1)\tilde\omega_1\wedge{\omega_0\over 1+g_k};$$
%in each case, we see that $\mathcal F$ is transversely affine.

%\noindent{\bf Subcase 2.2:} $\beta\neq0$. Therefore, one can write
%$\tilde\omega_N=h\beta$ for some meromorphic function $h\neq0$
%and we have 
%$$d\tilde\omega_N={dh\over h}\wedge\omega_N=0\ \ \ 
%\text{and thus}\ \ \ $$

\noindent{\bf Case 2.} Assume that $g_k\neq0$ for at least one 
$k\in\{0,2,\ldots,N-2\}$ but
$${1\over N-l}{dg_l\over g_l}={1\over N-k}{dg_k\over g_k}$$
for all $k,l\in\{0,2,\ldots,N-2\}$ such that $g_k,g_l\neq0$:
%For instance, when $N=3$, $g_k\neq0$ at most for $k=0$.
the closed $1$-form
$$\beta=\tilde\omega_1+{1\over N-k}{dg_k\over g_k}$$
does not depend on $k$.

\noindent{\bf Subcase 2.1:} $\beta=0$. Since 
$$\omega_0=\tilde\omega_1+g\cdot\tilde\omega_N,\ \ \ 
g=g_0+g_2+\cdots+g_{N-2}+1$$
we get that either $g=0$ and $\omega_0=\tilde\omega_1$ is closed, 
or $g\neq0$ and we have
$$d({\omega_0\over g})=d\tilde\omega_N=(N-1)\tilde\omega_1\wedge\tilde\omega_N
=(N-1)\tilde\omega_1\wedge{\omega_0\over g};$$
in each case, we see that $\mathcal F$ is transversely affine.

\noindent{\bf Subcase 2.2:} $\beta\neq0$. Therefore, one can write
$\tilde\omega_N=h\beta$ for some meromorphic function $h\neq0$
and we have 
$$d\tilde\omega_N={dh\over h}\wedge\tilde\omega_N.$$
Comparing with $d\tilde\omega_N=(N-1)\tilde\omega_1\wedge\tilde\omega_N$
and $\beta\wedge\tilde\omega_N=0$, we get
$$\left({dh\over h}-{N-1\over N-k}{dg_k\over g_k}\right)\wedge\tilde\omega_N=0.$$

\noindent{\bf Subsubcase 2.2.1:} ${N-1\over N-k}{dg_k\over g_k}={dh\over h}$
for all $k\in\{0,2,\ldots,N-2\}$ such that $g_k\neq0$. Then
$$\omega_0=\tilde\omega_1+gh\cdot\beta,\ \ \ 
g=g_0+g_2+\cdots+g_{N-2}+1$$
with $dg\wedge dh=0$. Since 
$\beta=\tilde\omega_1+{1\over N-1}{dh\over h}$, we get 
$$\omega_0=(1+gh)\tilde\omega_1 + {g\over N-1}dh.$$
Either $1+gh=0$ and $\omega_0$ is closed, or $1+gh\neq0$
and ${\omega_0\over 1+gh}$ is closed; in each case, 
$\mathcal F$ is transversely affine.

\noindent{\bf Subsubcase 2.2.2:} ${N-1\over N-k}{dg_k\over g_k}\neq{dh\over h}$
for at least one $k$. Therefore, we can conclude that
$$dg\wedge\tilde\omega_N=0\ \ \ \text{with}\ \ \ 
g:={h^{(N-k)}\over g_k^{(N-1)}}\ \ \ \text{non constant}$$
i.e. $\omega_N=fdg$ for some meromorphic function $f$. 
\end{proof}

\section{Examples}\label{S:Examples}

\subsection{Degree $2$ foliations on $\mathbb{C}P(n)$ have length $\le4$}

Here, we prove Theorem \ref{T:grau2}. In fact, given a degree $2$ 
foliation $\mathcal F$ on $\mathbb{C}P(n)$, we prove that,
after a convenient birational transformation
$$\Phi:\mathbb{C}P(n)\dashrightarrow\mathbb{C}P(n-1)\times\mathbb{C}P(1),$$ 
the tangency locus $\Delta$ between the foliation $\mathcal F':=\Phi_*\mathcal F$ 
and the projection $\pi:\mathbb{C}P(n-1)\times\mathbb{C}P(1)\to\mathbb{C}P(n-1)$
takes the following special form:
\begin{itemize}
\item either $\Delta$ is a vertical hypersurface, i.e. defined by $R\circ\pi=0$
for a non constant rational function $R$ on $\mathbb{C}P(n-1)$,
\item or $\Delta$ is the union of a vertical hypersurface like above
and the horizontal hyperplane at infinity $H_\infty:=\mathbb{C}P(n-1)\times\{\infty\}$.
\end{itemize}
One can easily deduce from this geometric picture that $\mathcal F'$
is actually defined by a unique rational integrable $1$-form 
$$\Omega=dz+\sum_{k=0}^N \omega_k z^k$$
where $\omega_k$ are rational $1$-forms on $\mathbb{C}P(n-1)$ 
and $z$ is the $\mathbb{C}P(1)$-variable. A Godbillon-Vey sequence 
of length $\le N+1$ is therefore provided by $(L_X^{(k)}\Omega)_k$
where $X=\partial_z$ is the vertical vector field. 
We will also prove that $N\le3$ in our case.
In the first case of the alternative above, we have $N\le2$:
$\Delta$ is vertical, $\mathcal F'$ is a Riccati foliation with respect
to $\pi$ and is in particular transversely projective. In the second case,
$N=2+m$ where $m$ is the multiplicity of contact between $\mathcal F'$
and the projection $\pi$ along the hyperplane at infinity $H_\infty$.
Actually, it is better to view $\Delta$ as a positive divisor, defined in charts
by the holomorphic function $\omega(X)$ where $X$ is a non vanishing holomorphic 
vector field tangent to the fibration given by $\pi$ and $\omega$ a holomorphic
$1$-form defining $\mathcal F'$ with codimension $\ge2$ zero set. Then, $m$ is the
weight of $\Delta$ along $H_\infty$.

Let $\mathcal F$ be a degree $2$ foliation on $\mathbb{C}P(n)$.
In order to construct $\Phi$ and reach the geometrical picture above, 
the rought idea is to find a rational pencil on $\mathbb{C}P(n)$ such that
the tangency locus $\Delta$ between the foliation and the pencil
intersects each rational fiber once. In fact, we choose any singular point $p$
of the foliation $\mathcal F$ and consider the pencil of lines passing through $p$.
Of course, the number of tangencies between a line and $\mathcal F$, counted
with multiplicities, is $2$, the degree of $\mathcal F$; but looking at the
pencil passing through $p$, we expect that the tangency occuring at the singular 
point disappear after blowing up the point $p$. Let us compute.

A foliation $\mathcal F$ of degree $\le2$ on $\mathbb{C}P(n)$ is given in an affine chart 
$\mathbb C^n\subset\mathbb{C}P(n)$ by a polynomial $1$-form with codimension $\ge2$
zero set having the special form
$$\Omega=\omega_0+\omega_1+\omega_2+\omega_3$$
where $\omega_i$ is homogeneous of degree $i$ and $\omega_3$ is radial (see
\cite{CerveauLinsNeto}):
we have $\omega_3(\mathcal R)=0$, 
where $\mathcal R:=x_1\partial_{x_1}+\cdots+x_n\partial_{x_n}$ is the radial vector
field. Saying that $\mathcal F$ is not of degree less than $2$ just means that,
if ever $\omega_3=0$, then $\omega_2$ is not radial.
Let us assume $p=0$ be singular for $\mathcal F$, i.e. $\omega_0=0$.
The tangency locus between $\mathcal F$ and the pencil of lines passing through $0$ 
is given by $\text{tang}(\mathcal F,\mathcal R)={\Omega(\mathcal R)=0}$.
If $\Omega(\mathcal R)$ is the zero polynomial, then this means that $\mathcal F$
is actually radial; we avoid this by choosing another singular point $p$.
Therefore, $\text{tang}(\mathcal F,\mathcal R)$ is a cubic hypersurface
which is singular at $p$. After blowing-up the origin, the foliation lifts-up
in the chart 
$$\pi:(t_1,\ldots,t_{n-1},z)\mapsto(zt_1,\ldots,zt_{n-1},z)=(x_1,\cdots,x_{n-1},x_n)$$
just by lifting-up the $1$-form $\Omega$ which now takes the special form
$$\pi^*\Omega=z\left((f_0(t)+zf_1(t))dz+z\tilde\omega_1+z^2\tilde\omega_2+z^3\tilde\omega_3\right)$$
where $f_0$ and $f_1$ are polynomial functions of $t=(t_1,\ldots,t_{n-1})$
and $\tilde\omega_i$ are polynomial $1$-forms depending only on $t$.
We observe that $\text{tang}(\mathcal F,\mathcal R)$ is now defined by
$\{z(f_0(t)+zf_1(t))=0\}$, has possibly some vertical components given
by common factors of $f_0$ and $f_1$ and has exactly $2$ non vertical components 
defined by $z=0$ and $z=-f_0/f_1$ (the two tangencies between any line of the pencil
with $\mathcal F$). Also, as expected, the first section $z=0$ is irrelevant 
since it disappears after division of $\pi^*\Omega$: the tangency locus between
the lifted foliation $\tilde{\mathcal F}$ and the lifted pencil 
(the vertical line bundle $\{t=\text{constant}\}$)
actually reduces to $\{f_0(t)+zf_1(t)=0\}$ in the chart above.
We now discuss on this set.

If $f_0\equiv0$, then ${\Omega\over z^2f_1(t)}$ is Riccati with wingular set over 
$\{f_1(t)=0\}$: $\mathcal F$
has length $\le3$. Recall that we have supposed $\mathcal F$ non radial 
and thus $f_0$ and $f_1$ cannot vanish identically simultaneously.

If $f_0\not\equiv0$, then the non vertical component of $\text{tang}(\mathcal F,\mathcal R)$
is the section $z=s(t)$, $s(t):=-{f_0(t)\over f_1(t)}$. If $f_1\equiv0$, then this section
is the hyperplane at infinity $\{z=\infty\}$: ${\Omega\over f_0(t)}$ is already
in the expected geometrical normal form and has length $\le4$. If $f_1\not\equiv0$, 
it suffices to push it towards infinity by a meromorphic change of coordinate of the form 
$\tilde z:={z\over z-s(t)}$; after this birational transformation, we are in the previous
case $y_1\equiv0$ and we have done. Precisely, the foliation is defined by
$$d\tilde z- \tilde z{\tilde\omega_1\over f_0}+\tilde z^2({2\tilde\omega_1\over f_0}-{df_0\over f_0f_1}
+{df_1\over f_1^2}-{\tilde\omega_2\over f_1})-\tilde z^3({\tilde\omega_1\over f_0}-{\tilde\omega_2\over f_1}
+{f_0\tilde\omega_3\over f_1^2}).$$

In order to finish the proof of Theorem \ref{T:grau2}, 
we note that a generic degree $2$ foliation of $\mathbb C P(2)$ 
has length $4$, i.e. is not transversely projective.
Actually, this is a well known fact. For instance,
it immediately follows from Corollary \ref{C:singularSimplyConnected}
and the fact that a generic degree $d\ge2$ foliation on $\mathbb C P(2)$
has no invariant algebraic curve. An explicit example is given in Section
\ref{S:Jouanolou}.

\begin{remark}\rm
If   $\mathcal F$ is a foliation of $\mathbb CP(2)$ given by a $1$-form of the type
$\omega = \omega_{\nu} + \omega_{\nu + 1} + f_{\nu + 1} (xdy - ydx)$ then, 
for generic $\omega$ as above, $\text{Tang}(\mathcal F,\mathcal E)$ is a rational curve 
and an argument similar to the one used above
implies that $\mathcal F$ also satisfies $l(\mathcal F)\le 3$.
\end{remark}

\eject

\subsection{The examples of Jouanolou}\label{S:Jouanolou}

In \cite{Jouanolou}, Jouanolou exhibited the first examples of holomorphic foliations of the projective
plane without algebraic invariant curves. His examples, one for each degree greater than or equal to $2$,
are the foliations of $\mathbb C P(2)$ induced by the homogeneous $1$-forms in $\mathbb C^3$
\[
 \Omega_n = \det \left( \begin{matrix} dx & dy & dz \\ x & y & z \\ y^n & z^n & x^n \end{matrix} \right) . 
\]
The automorphism group of the foliation $\mathcal J_d$, induced by $\Omega_d$, 
is isomorphic to  a 
semi-direct product of $\mathbb Z / (n^2 + n + 1)\mathbb Z$ with $\mathbb Z / 3 \mathbb Z$ 
and is generated 
by the transformations $\psi_n(x:y:z) = (\delta^{n^2} x :\delta^{n} y : \delta z )$ 
and $\rho(x:y:z)=(y:z:x)$, where $\delta$ is a primitive $(n^2+n+1)^{\text{th}}$ root 
of the unity.

In \cite{MaMoNoS} it is observed that the foliations $\mathcal J_n$ 
can be presented in a different way. 
If  $\mathcal F_n$ is the degree $2$ foliation of $\mathbb C P (2)$ induced by the $1$-form
\[
  \omega_n =   \det \left( \begin{matrix} dx & dy & dz \\ x & y & z \\ x(-x + ny) & y(-y + nz) &
  z(-z+nx) \end{matrix} \right) \, ,
\]
and $\phi_n:\mathbb C P(2) \dashrightarrow \mathbb C P(2)$ is the rational map
(of degree $n^2+n+1$) given by
\[
   \phi_n (x:y:z) = (y^{n+1}\cdot z:z^{n+1}\cdot x:x^{n+1}\cdot y) \, 
\]
then the foliation $\mathcal J_n$ is the pull-back of the foliation $\mathcal F_n$ 
under $\phi_n$, i.e., 
$\mathcal J_n = \phi_n^* \mathcal F_n$. Conversely  we can say that $\mathcal F_n$  
is birationally equivalent to the quotient
of $\mathcal J_n$ by the group generated by $\psi_n$.

From the results of the previous section it follows that  $\mathcal F_n$  has length at most
$4$. Pulling-back a 
Godbillon-Vey sequence by $\phi_n$ we obtain that the length of $\mathcal J_n$  is also
bounded by $4$ and since
it does not admit invariant algebraic curves its length is precisely $4$. 
We have therefore proved the

\begin{cor}
The foliations $\mathcal J_n$, for every $n\ge 2$, have length $4$.
\end{cor}

\eject

\subsection{A new component of the space of foliations on $\mathbb CP(3)$}\label{S:component}

We start by considering the transversely projective foliation on $\mathbb CP(2)$ 
given in the affine chart $\{(x,y)\} = \mathbb C^2 \subset \mathbb CP(2)$ 
by the $1$-form
\[
  \omega = x dy - ydx + P_2 dx + Q_2 dy + R_2(xdy - ydx) \, .
\]
where $P_2, Q_2, R_2$ are generic homogeneous polynomials of degree $2$.
This is a degree $2$ foliation of $\mathbb CP(2)$ 
transverse to the Hopf fibration $x/y = const$
outside three distinct lines.
Let us consider the homogenization $\Omega_3$ of $\omega$ 
in the coordinates $(x,y,z)$ of $\mathbb C^3$:
\[
  \Omega_3 = z^2 (x dy - ydx) + z(P_2 dx + Q_2 dy) + R_2(xdy - ydx) - R_3 dz \, ,
\]
where $R_3(x,y) = x P_2 + y Q_2$. 
The genericity condition on $P_2,Q_2,R_2$ implies that
 $d\Omega_3$ has only one zero
on $\mathbb C^3$ which is  isolated and located at the  origin. 
Of course, $\Omega_3$ defines a transversely projective foliation 
of $\mathbb C^3 \subset \mathbb CP(3)$. We will twist this foliation 
by a polynomial automorphism of $\mathbb C^3$. More precisely, 
if $\sigma(x,y,z)= (x,y,z+x^2)$ then
$$\Omega:=\sigma^* \Omega_3= \Omega_3 + \Omega_4 + \Omega_5\ \ \ \text{with}$$
$$\left\{\begin{matrix}
\Omega_3 &=& z^2(xdy - ydx) + z(P_2 dx + Q_2 dy) + R_2(xdy -ydx)- R_3 dz   \\
\Omega_4 &=& 2zx^2(xdy-ydx) + x^2(P_2dx+Q_2dy) - 2xR_3 dx\hfill \\
\Omega_5 &=&  x^4(xdy - ydx) \hfill
\end{matrix}\right.$$
%\begin{eqnarray*}
%\sigma^* \Omega_3 &=& (z+x^2)^2(xdy - ydx) + (z+x^2)(P_2 dx + Q_2 dy)  \\
%                  && +R_2(xdy - ydx) - R_3 dz -2R_3 xdx  \\
%                  &=& \{z^2(xdy - ydx) + z(P_2 dx + Q_2 dy) + R_2(xdy -ydx) 
%                 - R_3 dz \}  \\
%          && +\{2zx^2(xdy-ydx) + x^2(P_2dx+Q_2dy) - 2xR_3 dx \}  \\
%          && +\{ x^4(xdy - ydx) \} = \Omega_3 + \Omega_4 + \Omega_5 \,
%\end{eqnarray*}
The $1$-form $\Omega $ defines
a degree $4$ foliation on $\mathbb CP(3)$ which is transverse 
to the Hopf fibration(induced by the Euler vector field 
$x\frac{\partial}{\partial x} + y\frac{\partial}{\partial y} +
z\frac{\partial}{\partial z}$) outside the union of the four hyperplanes  
$\Omega_4(E)= x^2R_3(x,y)=0$. If $P_2,Q_2,R_2$
are generic, then these four hyperplanes are distinct.

Let $\mathcal F'$ be a foliation of degree $4$ close to $\mathcal F_\Omega$: 
$\mathcal F'$ is given in the affine chart
$\mathbb C^3$ by a polynomial $1$-form
\[
  \Omega' = \Omega_0' + \Omega_1' + \Omega_2' + \Omega_3' + \Omega_4' + \Omega_5' \, ,
\]
where the $\Omega_k'$ are homogeneous of degree $k$ and $\Omega_5'(E)\equiv 0$.

After normalization, we can suppose that the coefficients of $\Omega'$ 
are close to those of $\Omega$.
Since $d\Omega_3$ has an isolated singularity at $0$, 
there exists (see \cite{CamachoLinsNeto}) a point $0'$
where the $2$-jet of $\Omega'$ is zero, 
and the Euler vector field centered at $0'$ is in the kernel of the $3$-jet.
Therefore, after translating $0'$ to $0$, we can suppose that $\mathcal F'$ 
is given by
\[
  \Omega' = \Omega_3' + \Omega_4' + \Omega_5' \, .
\]
We verify that $\Omega'$ is transversely projective (with poles contained 
in $\Omega_4'(E)$).
In fact, since $\mathcal F$ is not transversely affine, 
the same holds for $\mathcal F'$.  
Therefore every element $\mathcal F'$ of the component of $\mathcal F(3,4)$
containing $\mathcal F$ is actually transversely projective.

\eject

\subsection{Transversely projective foliations that are not
pull-back}\label{S:modular}

\begin{example}[Example 8.6 of \cite{GhysGomez}]\rm
Let $\Gamma $ be discrete torsion free subgroup of $PSL(2,\mathbb R)^n$ such that the quotient
$ PSL(2,\mathbb R)^n/ \Gamma$ is compact. For $n\ge 2$, there exists examples 
such that the projection $\pi(\Gamma)$ on the first factor is a dense subgroup
of $PSL(2,\mathbb R)$ (see \cite{Borel}). 
The action of $\Gamma$ on $\mathbb H^n$, the $n$ product of the upper half-plane, 
is free, cocompact and preserves the regular foliation 
induced by the projection on the first factor.
In this way, we obtain a regular transversely projective foliation $\mathcal F$ 
on a $n$-dimensional compact complex manifold $M$ such that
every leaf is dense and the generic leaf is biholomorphic to $\mathbb H ^{n-1}$. 
Observe that  $\mathcal F$ is not the pull-back of a foliation 
on a lower dimensional manifold, 
otherwise there would exist compact subvarieties in $\mathbb H^{n-1}$. 
%Therefore, $\mathcal F$ satisfies the hypothesis of Theorem \ref{T:Hilbert}
%and its suspension examples on ruled manifolds of arbitrary dimension.
\end{example}

\begin{example}[Hilbert Modular Foliations]\rm
Let $K$ be a totally real number field of degree $n\ge2$ over the rational numbers $\mathbb Q$ 
and let $\mathcal O_K$ be the ring of integers of $K$. The group $\Gamma=PSL(2,\mathcal O_K)$
is dense in $PSL(2,\mathbb R)$, but considering the $n$ embeddings 
$i\circ\sigma:K\hookrightarrow\mathbb R$
given by the action $\sigma\in\text{Gal}(K/\mathbb Q)$,
we get an embedding $\Gamma\hookrightarrow PSL(2,\mathbb R)^n$
as a discrete subgroup of the product.
The quotient of $\mathbb H^n$, the $n$-product of the upper-half
plane $\mathbb H$, by $\Gamma$ is 
a quasiprojective variety $V$ which can be singular due to torsion elements of $\Gamma$. 
One can compactify and desingularize $V$ and obtain a projective manifold $M$.
The $n$ fibrations on $\mathbb H^n$ given by the projections on each of the factors
induce $n$ foliations on $M$ which are regular and pair-wise transversal 
outside the invariant hypersurfaces coming from the compactification and desingularization 
of $V$. By construction, they are transversely projective and all leaves apart
from the invariant hypersurface above are dense in $M$.
In \cite{Touzet} and \cite{MendesPereira}, 
some basic properties of these foliations are described.

When $K=\mathbb Q(\sqrt{5})$, the resulting variety is birationally equivalent 
to the projective plane. 
In \cite{MendesPereira} explicit equations for the foliations associated 
to the two projections $\mathbb H^2 \to \mathbb H$, denoted by
$\mathcal F_2$ and $\mathcal F_3$, are determined. 
We give below an explicit projective triple for them. 
The corresponding suspensions $\mathcal H_2$ and $\mathcal H_3$
defined by
$$\Omega=dz+\omega_0+z\omega_1+z^2\omega_2$$
can be seen as singular foliations on $\mathbb CP(2)\times\mathbb CP(1)$
or equivalently on $\mathbb CP(3)$. Although the leaves of $\mathcal F_2$ are dense,
we note that the same is not true for $\mathcal H_2$ since the monodromy lie in $PSL(2,\mathbb R)$.
\end{example}

\begin{landscape}
\begin{small}
{\bf A projective triple for $\mathcal F_2$}
\begin{eqnarray*}
\omega_0 &=&{\displaystyle \frac {(80\,y - 60\,x\,y - 80\,x^{2}
)dx + (36\,x^{2} - y - 32\,x) dy }{ - 720\,x^{3}\,y + 1728\,x^{5}
+ 80\,x\,y^{2} - y^{3} + 640\,x^{2}\,y - 1600\,x^{4}
 - 64\,y^{2}}}  \\
\omega_1&=&{\displaystyle \frac {16}{3}} \,{\displaystyle \frac
{24\,y + 3\,y^{2} + 62\,x\,y - 24\,x^{2} - 108\,x^{2}\,y -
368\,x^{3} + 432\,x^{4}}{ - 720\,x^{3}\,y + 1728\,x^{5} +
80\,x\,y^{2} - y^{3}
 + 640\,x^{2}\,y - 1600\,x^{4} - 64\,y^{2}}} dx \\
 &-& {\displaystyle \frac {4}{15}} \,{\displaystyle \frac {150\,y + 3
\,y^{2} + 192\,x - 172\,x\,y - 872\,x^{2} + 720\,x^{3}}{ - 720\,x
^{3}\,y + 1728\,x^{5} + 80\,x\,y^{2} - y^{3} + 640\,x^{2}\,y -
1600\,x^{4} - 64\,y^{2}}} dy \\
\omega_2 &=& {\displaystyle \frac {32}{45}} \frac{144\,y +
84\,y^{2} + 3\,y^{3}
 + 852\,x\,y - 48\,x\,y^{2} - 144\,x^{2} - 1472\,x^{2}\,y - 36\,x
^{2}\,y^{2} - 4416\,x^{3}
 + 336\,x^{3}\,y + 4928\,x^{4}}{
 - 720\,x^{3}\,y + 1728\,x^{5} + 80\,x\,y^{2} - y^{3} + 640\,x^{2
}\,y - 1600\,x^{4} - 64\,y^{2})} dx  \\
  &-& {\displaystyle \frac {32}{225}} \,{\displaystyle \frac {441\,y
+ 288\,x - 372\,x\,y + 12\,x\,y^{2} - 2292\,x^{2} - 56\,x^{2}\,y +
1472\,x^{3} - 108\,x^{3}\,y + 720\,x^{4}}{ - 720\,x^{3}\,y + 1728
\,x^{5} + 80\,x\,y^{2} - y^{3} + 640\,x^{2}\,y - 1600\,x^{4} - 64
\,y^{2}}} dy
\end{eqnarray*}
{\bf A projective triple for $\mathcal F_3$}
\begin{flushleft}
\begin{eqnarray*}
\omega_0 &=&  \frac { (- {\displaystyle \frac {5}{4}} \,y^{2} +
20\,x\,y - 60\,x^{3}) dx  +( - y + {\displaystyle \frac {3}{4}}
\,x\,y
 + x^{2}) dy}{ - 720\,x^{3}\,y + 1728\,x^{5} + 80\,x\,y^{2} - y^{3}
 + 640\,x^{2}\,y - 1600\,x^{4} - 64\,y^{2}} \\
\omega_1 &=& 4\,{\displaystyle \frac { - 80\,x^{2}\,y + 3\,y^{2} -
208\,x^{ 3} + 288\,x^{4} + 48\,x\,y}{ - 720\,x^{3}\,y +
1728\,x^{5} + 80\, x\,y^{2} - y^{3} + 640\,x^{2}\,y - 1600\,x^{4}
- 64\,y^{2}}} dx  \\
&-& {\displaystyle \frac {2}{5}} \,{\displaystyle \frac {40\,y +
y^{2 } - 168\,x^{2} + 192\,x^{3} - 50\,x\,y}{ - 720\,x^{3}\,y +
1728\, x^{5} + 80\,x\,y^{2} - y^{3} + 640\,x^{2}\,y - 1600\,x^{4}
- 64\, y^{2}}} dy \\
\omega_2 &=& {\displaystyle \frac {32}{5}} \,{\displaystyle \frac
{ - 9\,y ^{2} + 80\,x\,y + 8\,x\,y^{2} - 16\,x^{2}\,y - 368\,x^{3}
- 48\,x ^{3}\,y + 320\,x^{4}}{ - 720\,x^{3}\,y + 1728\,x^{5} +
80\,x\,y^{ 2} - y^{3} + 640\,x^{2}\,y - 1600\,x^{4} - 64\,y^{2}}}
dx \\
&+& {\displaystyle \frac {32}{25}} \,{\displaystyle \frac { - 36\,
y + 63\,x\,y + 164\,x^{2} - 28\,x^{2}\,y - 304\,x^{3} + 144\,x^{4
}}{ - 720\,x^{3}\,y + 1728\,x^{5} + 80\,x\,y^{2} - y^{3} + 640\,x
^{2}\,y - 1600\,x^{4} - 64\,y^{2}}} dy
\end{eqnarray*}
\end{flushleft}
\end{small}
\end{landscape}

\begin{thm}\label{T:Hilbert}
The explicit suspensions $\mathcal H_2$ and $\mathcal H_3$ above 
are not the meromorphic pull-back of a foliation on a surface.
\end{thm}

\begin{proof}
Suppose that there exists a foliation $\underline{\mathcal H}_2$ 
on a surface $S$
and a meromorphic map $\Phi:\mathbb CP(2) \times \mathbb CP(1)\dashrightarrow S$ 
such that $\Phi^* \underline{\mathcal H}_2 = \mathcal H_2$.

Let $U\subset \mathbb CP(2)\times\mathbb CP(1)$ be the Zariski open subset where $\Phi$ is holomorphic
and $U_0=U\cap(\mathbb CP(2) \times \{ 0 \})$.
After blowing-up $S$, one can assume $\Phi(U_0)$ having codimension $\le 1$.
The generic rank of $\Phi$ restricted to $U_0=U\cap(\mathbb CP(2) \times \{ 0 \})$
is $2$, otherwise we are in one of the following contradicting situations
\begin{enumerate}
\item  The closure of $\Phi(U_0)$ is a proper submanifold of $S$ non-invariant 
by $\underline{\mathcal H}_2$. In particular $\mathcal F_2$ is the pull-back of a foliation
of a foliation on a curve and is transversely euclidean; contradiction.
\item The closure of $\Phi(M_0)$ is a proper submanifold of $S$ invariant 
by $\underline{\mathcal H}_2$ (and not contained in the 
singular set
of $\underline{\mathcal H}_2$). Reasoning in local coordinates at the neighborhood
of a generic point $p\in\Phi(U_0)$, we see that
$\mathbb CP(2) \times \{ 0 \}$ is invariant by $\mathcal H_2$ obtaining a contradiction.
\end{enumerate}
We conclude therefore that $\Phi\vert_{U_0}$ is dominant and 
$\underline{\mathcal H}_2=\Phi_*\mathcal F_2$ has dense leaves 
(in fact all but finitely many). Therefore, the same density property holds
for the pull-back $\mathcal H_2=\Phi^*\underline{\mathcal H}_2$ providing a contradiction:
the Riccati foliation $\mathcal H_2$ has no dense leaf since its monodromy 
is contained in $PSL(2,\mathbb R)$. This proves the Theorem.
\end{proof}

\section{Integrable $1$-forms in  Positive Characteristic}

Due to the algebraic nature of many of the arguments used through this paper it is natural
to ask if it would be possible carry on a similar study for integrable $1$-forms on varieties
defined over fields of positive characteristic.

The surprising fact, at least for us, is that over fields of positive characteristic every $1$-form
admits a Godbillon-Vey sequence  of length one. In the case of $1$-forms on  the projective plane 
this is already implicitly proved in \cite{positiva}. 

Our argument is based on the following

\begin{lem}\label{L:pos}
Let $M$ be a $m$-dimensional smooth projective variety defined over an arbitrary field. If $\omega$
is an integrable rational $1$-form then there exists $m-1$ rationally independent vector fields $X_1, \ldots, X_{m-1}$ such that
\begin{enumerate}
\item $[X_i, X_j] = 0$ for every $i, j \in {1, \ldots, m-1}$;
\item $\omega(X_i) = 0$ for every $i \in {1, \ldots, m-1}$.
\end{enumerate}
\end{lem}
\begin{proof}
Let $f_1, \ldots, f_{m-1} \in k(M)$ be rational functions such that 
\[
\omega \wedge df_1 \wedge \cdots \wedge df_{m-1} \neq 0 \, .
\]
If $\omega_m=\omega$ and $\omega_i= df_i$, for $i=1 \ldots m-1$ then $\{\omega_i\}_{i=1}^{m}$ form 
a basis of the $k(M)$-vector space of rational $1$-forms over $M$. 

Let  $\{X_i\}_{i=1}^{m}$ be a basis of the space of rational vector fields on $M$ dual to $\{\omega_i\}_{i=1}^{m}$, i.e.,
$\omega_i(X_j) = \delta_{ij}$. It is clear that $\omega (X_i) = 0$ for every $i=1\ldots m-1$. We claim that $[X_i,X_j]=0$ for
every $i,j = 1 \ldots m-1$. It is  sufficient  to show that 
\begin{equation}\label{E:facil}
  \omega_k([X_i,X_j]) = 0 \, \, \text{for every} \, \, k=1\ldots m 
\end{equation}

For $k=m$ the integrability of $\omega$ implies that  (\ref{E:facil}) holds.
For $k<m$ we have that
\begin{eqnarray*}
  \omega_k ( [X_i,X_j]) &=& X_i(\omega_k(X_j)) - X_j(\omega_k(X_i)) + d\omega_k(X_i,X_j) =  \\
  &=& X_i(\delta_{kj}) - X_j(\delta_{ki}) + d^2 f_k(X_i,X_j) = 0 \, .
\end{eqnarray*}
This shows that (\ref{E:facil}) holds for every $k=1\ldots m $ and concludes the proof of the lemma.
\end{proof}

\begin{thm}
Let $M$ be a smooth projective variety defined over a field $K$ 
of characteristic $p > 0$
and $\omega$ be a rational $1$-form. If $\omega$ is integrable then $\omega$ admits
an "integrating factor", i.e., there exist a rational function $F \in K(M)$ such that
$F\omega$ is closed. Equivalently we have that 
\[
  d \omega = \omega \wedge \frac{dF}{F} \, .
\]
\end{thm}

\begin{proof}
Let $m$ be the dimension of $M$ and $X_1, \ldots, X_{m-1}$ be the rational vector fields 
given by lemma \ref{L:pos}. We will distinguish two cases:
\begin{enumerate}
\item for every $i=1\ldots m-1$ we have that $\omega(X_i ^p)=0$
\item there exists $i \in \{1,\ldots, m-1\}$ such that $\omega(X_i ^p)\neq 0$
\end{enumerate}

Let $\mathcal F$ be the unique saturated subsheaf of the tangent sheaf of $M$ which 
coincides with the kernel of $\omega$  over the generic point of $M$. The integrability of 
$\omega$ implies that $\mathcal F$ is involutive. If we are in the case (1) then we have
also that $\mathcal F$ is $p$-closed. From \cite[propositions 1.7 and 1.9, p. 55--56]{MiPe} it follows that $\omega=gdf$ where
$g, f \in k(M)$.

In case (2) we can suppose that $\omega(X_1 ^p)\neq 0$. If  $F= \omega(X_1 ^p)^{-1}$ then 
\[
   d ( F\omega) = F\omega \wedge L_{X_1 ^p} ( F \omega) \, .
\]

To conclude we have just to prove that $L_{X_1 ^p} ( F \omega) = 0$. 
In fact since $F \omega(X_1 ^p)=1$ it follows  that 
\[ 
   L_{X_1 ^p} ( F \omega) = i_{X_1^ p} d (F\omega) \, .
\]
Moreover for every $i=1\ldots m-1$ we have that $[X_1^p, X_i] =0$, since $X_1$ commutes with $X_i$, and therefore
\[
   i_{X_1^ p} d (F\omega)(X_i) = F \omega([X_1^p, X_i])  - X_1^p( F \omega(X_i) ) + X_i(F \omega(X_1^p) ) = 0 \, .
\] 
This is sufficient to show that $L_{X_1 ^p} ( F \omega) = 0$ concluding the proof of the Theorem.
\end{proof}

As a corollary we obtain a codimension one version of the main result of \cite{positiva}.

\begin{cor}
Let $\omega$ be a polynomial integrable $1$-form on $\mathbb A^n _k$, where $k$ is a field of positive characteristic. If $d\omega \neq 0$ then 
there exists an irreducible algebraic hypersurface H such that $i^* \omega=0$, where $i:H \to \mathbb A^n _k$ denotes the inclusion.
\end{cor}
\begin{proof}
Of course $\omega$ can be interpreted as rational $1$-form over $\mathbb P^n _k$ which is regular over 
$\mathbb A^n _k$. From Theorem \ref{T:positiva} there exists a rational function $F \in k(x_1,\ldots, x_n)$ such that 
\[
 d \omega = \omega \wedge \frac{dF}{F} \, .
\]
Since $d\omega \neq = 0$ we have that $dF \neq 0$, i.e., $F$ is not a $p$-th power. In particular
the polar set of $dF/F$ is not empty. It is an easy exercise to show that every irreducible component $H$
of the polar set of $dF/F$ satisfies $i^* \omega = 0$, where $i:H\to \mathbb A^n _k$ denotes the inclusion 
\end{proof}

In fact the same proof as  above yields the stronger

\begin{cor}
Let $\omega$ be a regular integrable $1$-form over a smooth quasiprojective algebraic variety $M$ defined over $k$, a  field of positive characteristic.
Suppose that  $ H^0(M,\mathcal O_M ^*) = k^*$. If  $d\omega \neq 0$ then 
there exists an irreducible algebraic hypersurface H such that $i^* \omega=0$, where $i:H \to M$ denotes the inclusion.
\end{cor}

Observe that the result above can be applied to projective varieties since there exists such varieties with global regular $1$-forms
which are not closed, see \cite{MiPe}.

\end{document}